\documentclass[11pt]{article}
\usepackage[pagewise]{lineno}
\usepackage[pdftex, dvips]{graphicx} 
\usepackage[utf8]{inputenc}
\usepackage{subfig}
\usepackage{epstopdf}
\usepackage{float}
\usepackage[english]{babel}
\usepackage{latexsym,amsmath,amssymb,amsfonts,cancel,dsfont}
\usepackage{bbm}
\usepackage{fancyhdr}
\usepackage[margin=2.5cm,top=3cm]{geometry}

\usepackage{cite}
\usepackage{csquotes}
\usepackage{xcolor}
\usepackage{caption}
\usepackage[colorlinks=true]{hyperref}

\newtheorem{defi}{Definition}[section]
\newtheorem{teo}[defi]{Theorem}
\newtheorem{prop}[defi]{Proposition}
\newtheorem{lem}[defi]{Lemma}

\newtheorem{ejem}[defi]{Example}
\newtheorem{remark}[defi]{Remark}

\newcommand{\dis}{\displaystyle}
\newcommand{\N}{\mathbb{N}}
\newcommand{\Q}{\mathbb{Q}}
\newcommand{\R}{\mathbb{R}}
\newcommand{\ep}{\varepsilon}

\newcommand{\la}{\lambda}

\newcommand{\va}{\varphi}

\newcommand{\black}{\color{black}}

\begin{document}

\title{ \huge 
\textbf{Uniqueness and numerical reconstruction\\
for inverse problems dealing with\\
interval size search}}

\author{
J. Apraiz\thanks{Universidad del Pa\'is Vasco, Facultad de Ciencia y Tecnolog\'ia, Dpto.\ Matem\'aticas, Barrio Sarriena s/n 48940 Leioa (Bizkaia), Spain. E-mail: {\tt jone.apraiz@ehu.eus}.},
\ \ J. Cheng\thanks{Fudan University, Shanghai 200433, China,  E-mail: {\tt jcheng@fudan.edu.cn}.},
\ \ A. Doubova\thanks{Universidad de Sevilla, Dpto.\ EDAN e IMUS, Campus Reina Mercedes, 41012~Sevilla, Spain, E-mail: {\tt doubova@us.es}.},
\ \ E. Fern\'andez-Cara\thanks{Universidad de Sevilla, Dpto.\ EDAN e IMUS, Campus Reina Mercedes, 41012~Sevilla, Spain, E-mail: {\tt cara@us.es}.}, 
\ \ M. Yamamoto\thanks{University of Tokyo, Japan, E-mail: {\tt myama@next.odn.ne.jp}.}
}

\date{}

\maketitle

\begin{abstract}
We consider a heat equation and a wave equation
in a spatial interval over a time interval.  
This article deals with inverse problems of determining sizes of spatial 
intervals by extra boundary data of solutions of the governing equations.
Under several different circumstances, we prove the uniqueness, 
the non-uniqueness and some size estimate.
Moreover, we numerically solve the inverse problems and compute accurate 
approximations of the sizes.  This is illustrated with satisfactory numerical 
experiments.
\end{abstract}

\vspace*{0,5in}
 
\textbf{AMS Classifications:} {\black 35R30, 35K05, 35L05, 65M32.}

\textbf{Keywords:} {\black  Inverse problems, uniqueness, heat equation, wave equation, numerical reconstruction.}

\section{Introduction}\label{Sec-1}

   In this article, we study inverse problems where the goal is to determine or estimate the size of the spatial interval where the governing equation holds.
   Although we can discuss more general cases, we concentrate here on the one-dimensional heat and wave equations.
   
   In recent years, the interest for the analysis and solution of inverse problems of many kinds has grown a lot.
   This is motivated for their relevance in many applications: elastography and medical imaging, seismology, potential theory, ion transport problems or chromatography, finances, etc.;
   {\black see for instance~\cite{Borcea_et_al, Hanke, Richter}.}

   In order to understand the situation, let us begin by recalling what is, in our framework, a {\it direct} problem.
   In general terms, in a direct problem we try to find (exactly or approximately) one or several functions that model a phenomenon.
   In these problems, the geometrical data, the media properties (expressed by the coefficients in the equations) and the initial and boundary data are assumed to be known;
   we use them to compute the solution to a governing system and, then, we obtain some generally useful information.
   
   In a related inverse problem, only some data (and not all of them) are known but we have access to this additional information (the observation).
   Thus, the final goal is to recover or estimate the unknown data, which will make it possible to compute the solution.
   Actually, the variety of inverse problems is huge compared with direct problems and several kinds of inverse problems which come from very classical and basic direct problems, still wait for theoretical and numerical researches.
Let us mention for example monographs
Bellassoued and~Yamamoto~\cite{BelYam3},
~Isakov~\cite{Isakov},
Romanov~\cite{Romanov} and~Hasanov and~Romanov~\cite{HasRom},
where various theoretical and numerical aspects of inverse problems 
for partial differential equations are described.

Among various types of inverse problems, our article focuses on inverse 
problems of determining sizes of spatial intervals where the heat equations and
wave equations hold.
More precisely, we consider 1D heat and wave equations in 
~$(x,t) \in (0,\ell)\times (0,T)$.
In the inverse problems considered in this article, 
we try to determine the width $\ell$ of the spatial $x$-interval from extra 
data.
This is much less studied, although it is connected to many interesting 
applications, such as non-destructive testing.
As for this kind of inverse problems related to the heat equation,   
we refer for example to 
Banks, Kojima and~Winfree~\cite{Banksk}, Bryan and Caudill~\cite{Bry1}--\cite{Bry3}, 
Chapko, Kress and~Yoon~\cite{Chap1, Chap2}, Fredman~\cite{Fred}, 
Wang, Cheng, Nakagawa and~Yamamoto~\cite{WaChNaY} 
and~Wei and~Yamamoto~\cite{WeiY}. 

   Very frequently, inverse problems are not {\it well-posed} in the sense of~\cite{Hadamard}.
   This means that, either the solution does not exist, or it is not unique and/or small errors in the observation or known data lead to large errors in the solution(s).
   In particular, uniqueness is the first theoretical issue for our inverse problem.  
   As our main results show, the uniqueness and the non-uniqueness are sensitive according to data settings.
   The main techniques used here are usual in control and parameter  identification theory:
   tools from real and complex analysis, Fourier expansions, known estimates (energy or Carleman-like), multipliers method, etc.                                    

The second main issue of this article is the reconstruction method for our 
inverse problem.
   Here, the goal is to compute approximations of lengths as solutions to the inverse problems using the observation data.
   To our best knowledge, there are no works combining the theoretical studies on the uniqueness and the numerical reconstruction of sizes of intervals for a heat and a wave equation.

   As shown below, an efficient technique to achieve this relies on  a reformulation of the search of the length as an extremal problem.
   This is classical nowadays and has been applied in a lot of situations; {\black see for instance~Lavrentiev and~others~\cite{Lavrentiev_et_al}, Samarskii and~Vabishchevich~\cite{Samarskii} and~Vogel~\cite{Vogel}.}
   
\
   
The organization of this article is as follows.
   
   In Section~\ref{Sec-2}, we consider the 1D heat equation.
   From the initial and boundary data and the outgoing heat flow, we try to 
   recover the length of the spatial interval where the problem is posed.
   We separately consider two cases: zero initial data and non-zero
   initial data.  In the case of zero initial data, a classical 
   argument yields uniqueness on the basis of 
   the unique continuation property.
   On the other hand, if initial data are not zero, then the uniqueness 
   is more delicate.  More precisely, in ~Section~\ref{SSec-2.1}, 
   in the case where initial data are not zero and Dirichlet boundary input is 
   zero, this will be completed  with some properties satisfied by 
   lenghts leading to the same observation, and we will provide 
   necessary conditions for the non-uniqueness.
   Then, in~Section~\ref{SSec-2.2}, we will analyze the case where both 
   boundary and initial data are non-zero.
   We will see that the size of the boundary data is determinant for a positive    uniqueness result: our result tells that 
   large boundary input can guarantee uniqueness.

   In Section~\ref{Sec-3}, we will deduce results for the wave equation.
   Many of them are like in the case of the heat equation, 
   although the energy conservation forces some necessary changes.
   In particular, for large boundary data we only get now asymptotic uniqueness.
   
   Section~\ref{Sec-4} deals with the results of several numerical experiments.
   They have been implemented to illustrate the theoretical results in the 
   previous sections.
   The solutions to the inverse problems have been computed in a nowadays 
   standard way, by introducing appropriate extremal problem reformulations. In the numerical tests, we have used 
   a MatLab Optimization Toolbox function as a solver.
   
   Finally, in~Section~\ref{Sec-5}, we have collected some additional comments 
   and questions.
   
\

   Throughout this paper, $\|\cdot\|$ and~$(\cdot\,,\cdot)$ will stand for the  usual $L^2$ norm and scalar product, respectively.
   In the particular case of the space $L^2(0,\ell)$, we will sometimes write $(\cdot\,,\cdot)_{\ell}$ in order to make explicit the length~$\ell$.
   The symbol $C$ will denote a generic positive constant.

\section{Some positive and negative results for the heat equation}\label{Sec-2}

   We begin with considering the 1D classical heat equation in a rod of length $\ell$ with Dirichlet data:
   \begin{linenomath*}
\begin{equation}\label{problema}
\begin{cases}
u_t - u_{xx} = 0,\, & 0 < x < \ell,\, 0 < t < T,\\
u(0,t)=\eta(t),\ \ u(\ell,t)=0,\, & 0 < t < T,\\
u(x,0)=u_0(x),\, & 0<x<\ell.
\end{cases}
\end{equation}
   \end{linenomath*}

   \textbf{Problem~IP-1 (heat equation, observation on the known boundary):}
\textit{Fix $u_0 = u_0(x)$ and $\eta = \eta (t)$ in appropriate spaces and assume that $u_x|_{x=0}$ (the heat flow on the left) is known. Then, find $\ell$.}

\

   We interpret $u_0$ and~$\eta$ as inputs for our inverse problem, while the~$u_x(0,t)$ with~$0<t<T$ are the observed results.

   We are interested in proving uniqueness.
   More precisely, the following question is in order:
  
\

\textbf{Uniqueness (observation on the known boundary):} \textit{Let $u^{\ell}$ and $u^L$ be the solutions to~\eqref{problema} corresponding to the spatial intervals $(0,\ell)$ and $(0,L)$, respectively. Assume that the corresponding observations $u_x^\ell(0,\cdot)$ and~$u_x^L(0,\cdot)$ coincide, that is,
   \[
u^{\ell}_x(0,t) = u^{L}_x(0,t) \ \text{ in } \ (0,T).
   \]
   Then, do we have $\ell=L$?}

\

   For $u_0 \in L^2(0,L)$, we note that the series of the $x$-derivatives at zero of the terms in the expansions \eqref{sol1} and~\eqref{sol2} converge in~$C^0([\epsilon,T])$ for arbitrary $\epsilon>0$.

\

In the sequel, we will provide some positive and negative answers to this question, depending on the kind of imposed boundary data $\eta$ or initial data $u_0$. We will also present some results concerning stability estimates.


\subsection{Zero initial and/or boundary data}\label{SSec-2.1}

\subsubsection{ \textbf{Case I}: $\eta\not\equiv 0$ and $u_0\equiv 0$}\label{SSSec-2.1.1}

In this case, the uniqueness is a well-known result by the unique continuation 
(e.g., Isakov \cite{Isakov}, Vessella \cite{V}, Yamamoto \cite{Ya}).  
Since the argument is well-known, we only sketch the essence. 

\begin{itemize}

\item Let $\ell$ and~$L$ be given with $0 < \ell < L$ and let us assume that the associated solutions to~\eqref{problema}, $u^\ell$ and~$u^L$, satisfy $u_x^\ell(0,\cdot) = u_x^L(0,\cdot)$.

\item Let us introduce $v := u^\ell - u^L$. Then $v$ solves the heat equation in~$(0,\ell) \times (0,T)$, $v(0,t) = 0$ and~$v_x(0,t) = 0$ in~$(0,T)$. Therefore, from the unique continuation property, one has $v \equiv 0$ in~$(0,\ell) \times (0,T)$.

\item This implies in particular that $u^L(\ell,t) \equiv 0$ and, since $u_0 = 0$, $u^L = 0$ in~$(\ell,L) \times (0,T)$. Using again unique continuation, we deduce that $u^L \equiv 0$. But this is an absurd, since we had $\eta\not\equiv 0$.

\end{itemize}

\subsubsection{ \textbf{Case II}: $\eta\equiv 0$ and $u_0\not\equiv 0$}\label{SSSec-2.1.2} 

   In general, when the left boundary condition is zero and we impose nonzero initial data, we can prove that uniqueness does not hold. For this purpose, we will use the eigenvalues and eigenfunctions of the Dirichlet-Laplace operators in~$(0,\ell)$ and~$(0,L)$:
   \begin{linenomath*}
\begin{equation*}\label{eigenv}
\left\{ \begin{array}{l} \dis
\la_n := \frac{n^2\pi^2}{\ell^2}, \quad n \in \N, 
\\ \dis
\va_n(x) :=\sqrt{\frac{2}{\ell}}\sin \frac{n\pi x}{\ell}, \quad n \in \N, \quad 0<x<\ell
\end{array}\right.
\end{equation*}
   \end{linenomath*}
and
   \begin{linenomath*}
\begin{equation*}\label{eigenf}
\left\{ \begin{array}{l} \dis
\mu_n := \frac{n^2\pi^2}{L^2}, \quad n \in \N,
\\ \dis
\psi_n(x) := \sqrt{\frac{2}{L}}\sin \frac{n\pi x}{L}, \quad n \in \N, \quad 0<x<L. 
\end{array}\right.
\end{equation*}
   \end{linenomath*}
   The solutions to~\eqref{problema} with $\eta\equiv 0$ corresponding to~$\ell$ and~$L$ can be defined for all~$t > 0$. They are respectively given by
   \begin{linenomath*}
\begin{align}
u^{\ell}(x,t)=& \sum_{n=1}^{\infty}  (u_0, \va_n)_{\ell} \va_n(x) e^{-\la_nt},\quad 0<x<\ell,\, t>0\label{sol1}
\\ \noalign{\noindent\mbox{and}}
u^{L}(x,t)=& \sum_{n=1}^{\infty} (u_0, \psi_n)_{L} \psi_n(x) e^{-\mu_nt},\quad 0<x<L,\, t>0,\label{sol2}
\end{align}
   \end{linenomath*}
where we have introduced the scalar products
   \begin{linenomath*}
\begin{equation*}
(f,g)_{\ell} := \int^{\ell}_0 f(x)g(x) \,dx \quad\text{and}\quad
(f,g)_L := \int^L_0 f(x)g(x) \,dx.
\end{equation*}
   \end{linenomath*}

Henceforth, by $\#{ \{n : (u_0, \va_n)_{\ell} \ne 0\} }$ we denote the 
number of the elements of the set.

The following holds:

\begin{prop}\label{prop2.1}
   If $L/\ell \in \Q$, then there exist initial data $u_0$ verifying 
   \begin{linenomath*}
\begin{equation}\label{con1}
\#{ \{n : (u_0, \va_n)_{\ell} \ne 0\} } = \#{ \{n : (u_0, \psi_n)_{L} \ne 0\} } =1,
\end{equation}
   \end{linenomath*}
such that $u^{\ell}_x(0,t) = u^{L}_x(0,t)$ for all~$t>0$. Thus, we can have non-uniqueness with initial data $u_0$ satisfying~\eqref{con1} 
even if $|L-\ell|$ is arbitrarily small.
\end{prop}

\textbf{Proof}:
   Let $m_0, n_0 \in \N$ be given, such that $n_0 < m_0$ and $\ell = n_0 L/m_0$, that is, $m_0/L = n_0/\ell$. Let us choose $k_1, n_1 \in \N$ such that $n_1 = k_1 m_0/n_0$.
   Note that
   \begin{linenomath*}
\begin{equation*}
\la_{k_1} = \frac{k_1^2\pi^2}{\ell^2} = \frac{n_1^2\pi^2}{L^2}
= \mu_{n_1}
\end{equation*}
   \end{linenomath*}
and set
   \begin{linenomath*}
\begin{equation*}
u_0(x) := \sin \frac{k_1\pi x}{\ell} = \sin \frac{n_1\pi x}{L}, \quad x\in \R.
\end{equation*}
   \end{linenomath*}

   The functions in~\eqref{sol1} and~\eqref{sol2} corresponding to this $u_0$ are the following:
   \begin{linenomath*}
\begin{equation*}
u^{\ell}(x,t) = \sqrt{\frac{2}{\ell}} (u_0,\va_{k_1})_{\ell}\sin \frac{k_1\pi x}{\ell}e^{-\la_{k_1}t} 
= \sin \frac{k_1\pi x}{\ell}e^{-\la_{k_1}t}
\end{equation*}
   \end{linenomath*}
and
   \begin{linenomath*}
\begin{equation*}
u^{L}(x,t) = \sqrt{\frac{2}{L}} (u_0,\psi_{n_1})_{\ell}\sin \frac{n_1\pi x}{L}e^{-\mu_{n_1}t}
= \sin \frac{n_1\pi x}{L}e^{-\mu_{n_1}t}.
\end{equation*}
   \end{linenomath*}
   Consequently,
   \begin{linenomath*}
\begin{equation*}
u^{\ell}_x(x,t)=\frac{k_1\pi}{\ell}\cos \frac{k_1\pi x}{\ell}e^{-\la_{k_1}t} \ \text{ and } \ 
u^{L}_x(x,t)=\frac{n_1\pi}{L}\cos \frac{n_1\pi x}{L}e^{-\mu_{n_1}t}
\end{equation*}
   \end{linenomath*}
and we conclude that $u^{\ell}_x(0,t) \equiv u^{L}_x(0,t)$. \hfill $\square$

\

   Next, we will show that, if the solutions $u^\ell$ and $u^L$ produce the same observation on the left boundary, i.e.~$u^{\ell}_x(0,\cdot) 
\equiv u^{L}_x(0,\cdot)$, then the associated lengths~$\ell$ and~$L$ 
cannot be completely independent:

\begin{prop}\label{prop2.2}
   Assume that $u^{\ell}_x(0,t) = u^{L}_x(0,t)$ in~$(0,T)$ and set $n_0 :=  \min\{ n : (u_0, \va_n)_{\ell} \ne 0\} \geq 1$. Then, there exists $N\in \N$ such that
   \begin{linenomath}
\begin{equation}\label{6p}
L = \frac{N}{n_0}\ell.\\
\end{equation}
   \end{linenomath}
   In particular, if $n_0=1$, one has $L = N\ell$, which means that only 
multiples of $\ell$ can provide the same observation at~$x=0$.
   Furthermore, if we assume that $\vert L-\ell\vert < \delta$ and  $L, \ell > \delta_0$, where  $\delta, \delta_0 >0$ and~$\delta < \delta_0/n_0$, then
   \[
u_x^\ell(0,t) = u_x^L(0,t) \ \text{ for } \ 0<t<T \ \Rightarrow \ L=\ell.
   \]
\end{prop}

\

The final part of the proposition means that if we a priori know that $L-\ell$ is sufficiently small compared with $\ell$ and~$L$, then uniqueness holds.

\

\textbf{Proof}:
   From~\eqref{sol1} and~\eqref{sol2}, we see that
   \begin{linenomath*}
\begin{equation*}
\frac{1}{\ell^{3/2}} \sum_{n=1}^{\infty} n (u_0,\va_n)_{\ell} e^{-\la_nt}
= \frac{1}{L^{3/2}} \sum_{n=1}^{\infty} n (u_0,\psi_n)_L e^{-\mu_nt} .
\end{equation*}
   \end{linenomath*}
   Taking Laplace transforms in both sides and applying analytic continuation, we find that
   \begin{linenomath*}
\begin{equation}\label{equlap}
\sum_{n=1}^{\infty} \frac{n(u_0,\va_n)_{\ell}}{\ell^{3/2}} 
\frac{1}{z+\la_n} 
= \sum_{n=1}^{\infty} \frac{n(u_0,\psi_n)_L}{L^{3/2}} \frac{1}{z+\mu_n}
\ \mbox{ for $z \not\in \{-\la_n\}_{n\in \N} \cup \{-\mu_n\}_{n\in \N}$.}
\end{equation}
   \end{linenomath*}
   
   Let us assume that $-\la_{n_0} \not\in \{-\mu_n\}_{n\in \N}$. Then, integrating both sides of~\eqref{equlap} on a small circle centered at $-\la_{n_0}$ and applying Cauchy's Theorem, we obtain $(u_0,\va_{n_0})_{\ell} = 0$, which is impossible by the definition of $n_0$. Therefore, there exists $N \in \N$ such that 
   \begin{linenomath*}
\begin{equation*}\label{lamu}
\la_{n_0} = \mu_N,                    
\end{equation*}
   \end{linenomath*}
that is, $n_0/\ell = N/L$ and~\eqref{6p} holds. 

We prove the final part as follows.  By (5), we have 
$$
\left\vert \frac{N-n_0}{n_0}\ell\right\vert 
= \vert L-\ell\vert
< \delta < \frac{\delta_0}{n_0} < \frac{\ell}{n_0}.
$$
Therefore, $\vert N-n_0\vert < 1$.  Since $N-n_0$ is an integer, this means
that $N-n_0=0$, that is, $L=\ell$.
\hfill $\square$

\

   In the following result, we give a characterization of~$L$ in terms of~$\ell$ assuming that $(u_0, \va_n)_{\ell} \ne 0$ for all $n$:

\begin{prop}\label{prop2.3}
Assume that $u^{\ell}_x(0,t) = u^{L}_x(0,t)$ in~$(0,T)$ and $(u_0, \va_n)_{\ell} \ne 0$ for all $n \geq 1$. Then, there exists $N\in \N$ such that 
\begin{linenomath*}
\begin{equation*}
L = N\ell.
\end{equation*}
\end{linenomath*}
\end{prop}

   We see that the conclusion is stronger than in Proposition~\ref{prop2.2}. Before proving Proposition~\ref{prop2.3}, we give an example that explains the effects of the symmetry of initial data.

\begin{ejem}
   Assume that $\ell > 0$ and $L=2\ell$.  In this case, if
\begin{linenomath*}
\begin{equation*}
u_0(x) = u_0(L-x), \quad 0<x<L,
\end{equation*}
\end{linenomath*}
one has $u^{\ell}_x(0,t) \equiv u^L_x(0,t)$, that is, the same observation is obtained at $x=0$.
\end{ejem}

\textbf{Proof of Proposition~\ref{prop2.3}}:
   Arguing as in the proof of Proposition~\ref{prop2.2}, we see that, for each  $n \in \N$, there exists $m(n) \in \N$ such that 
$\la_n = \mu_{m(n)}$, that is,
   \begin{linenomath*}
\begin{equation*}
\frac{L}{\ell} = \frac{m(n)}{n} \quad \mbox{for all $n \in \N$.}
\end{equation*}
   \end{linenomath*}
   Hence, $m(n)/n$ is independent of $n$ and we can introduce $r := m(n)/n = m_0/n_0$, where $m_0, n_0 \in N$ have 
no common divisor other than $1$.
   Observe that $r\in \Q$ and $rn \in \N$ for all $n\in \N$.

   Assume that $n_0 \ne 1$.  We arbitrarily choose $\widetilde{n} \in \N$
which has no common divisors other than $1$ with $n_0$.
   Then $r\widetilde{n} = \frac{m_0\widetilde{n}}{n_0} \in \N$.
   This is impossible by $n_0\ne 1$.

   Therefore, $r \in \N$ and, we finally have $L = N\ell$ with $N = r$, which completes the proof.\hfill $\square$

\subsection{Results where $\eta (t)\not\equiv 0$ and $u_0(x)\not\equiv 0$}\label{SSec-2.2}

   The situation is much more complex if we allow nonzero $\eta$ and $u_0$. Actually, we will see in this section that, if $\eta$ is ``large enough'', uniqueness hold.

   To this purpose, some auxiliary lemmas are needed. In the first one, we recall a conditional stability property of the Cauchy problem for the 1D heat equation:

{\black
\begin{lem}\label{lemma2.8}
   Let $L_*$, $T_0$, $T_1$ and~$\varepsilon_0$ be positive constants and assume that $0 < T_0 < T_1$ and~$0 < \varepsilon_0 \leq \varepsilon < L/2 \leq L_*/2$. Let~$v$ be a solution to the heat equation in~$(0,L) \times (T_0,T_1)$, with
   \[
\| v(L,\cdot)\|_{L^2(T_0,T_1)}  + \| v_x(L,\cdot)\|_{L^2(T_0,T_1)} \le D
   \]
and
   \[
\| v_x\|_{L^{\infty}((\ep_0,L_*)\times (T_0,T_1))} \le M.
   \]
   There exist constants $K>0$ and $\theta \in (0,1)$ (depending only 
on~$\varepsilon_0$, $L_*$, $T_0$, $T_1$ and~$M$) such that 
\begin{linenomath*}
\begin{equation*}
\left( \int^{T_1}_{T_0} \int^{2\varepsilon}_{\varepsilon} 
\| v(x,t)\|^2 dxdt \right)^{\frac{1}{2}} \le K D^{\theta}.
\end{equation*}
   \end{linenomath*}
\end{lem}

   The proof follows classical arguments relying on Carleman estimates and can be found in~\cite{Ya}.
}

   The second auxiliary result concerns traces of functions in~$H^2(0,\ell)$:

\begin{lem}\label{traceteo}
   Let $L_* > 0$ be given.
   Then
   \begin{linenomath*}
\begin{equation*}
\left| \frac{df}{dx}(0) \right| \le \frac{C(L_*)}{\ell^{3/2}}
\,\| f\|_{H^2(0,\ell)}
\end{equation*}
   \end{linenomath*} 
for any $f \in H^2(0,\ell)$ and any~$\ell \in (0,L_*)$.
\end{lem}

\noindent
\textbf{Proof:}
   Let us introduce $y := x/\ell$ and $g(y) := f(x)$. Then, $x\in (0,\ell)$ if and only if $y \in (0,1)$ and the following identities hold:
   \begin{linenomath*}
\begin{equation*}
\int^{\ell}_0 |f(x)|^2 \,dx = \ell\int^1_0 |g(y)|^2 \,dy,
\quad 
\int^{\ell}_0 \left| \frac{df}{dx}(x)\right|^2 \,dx 
= \frac{1}{\ell}\int^1_0 \left| \frac{dg}{dy}(y)\right|^2 \,dy
\end{equation*}
   \end{linenomath*} 
and
   \begin{linenomath*}
\begin{equation*}
\int^{\ell}_0 \left| \frac{d^2f}{dx^2}(x)\right|^2 \,dx 
= \frac{1}{\ell^3}\int^1_0 \left| \frac{d^2g}{dy^2}(y)\right|^2 \,dy.
\end{equation*}
   \end{linenomath*} 
   Consequently, the Trace Theorem in the spatial interval $(0,1)$ implies the existence of a pure constant $C_0 > 0$ such that
   \begin{linenomath*}
\begin{equation*}
\ell^2 \left| \frac{df}{dx}(0)\right|^2 
= \left| \frac{dg}{dy}(0)\right|^2 
\le C_0 \int^1_0 \left( |g(y)|^2 + \left| \frac{dg}{dy}(y)
\right|^2 + \left| \frac{d^2g}{dy^2}(y)\right|^2 \right) \,dy,
\end{equation*}
   \end{linenomath*} 
which can also be written in the form
   \begin{linenomath*}
\postdisplaypenalty=0
\begin{align*}
\left| \frac{df}{dx}(0)\right|^2 
&\le \frac{C_0}{\ell^2}
\int^{\ell}_0 \left( \frac{1}{\ell} |f(x)|^2
+ \ell \left| \frac{df}{dx}(x)\right|^2
+ \ell^3 \left| \frac{d^2f}{dx^2}(x)\right|^2 \right) \,dx \\
&\le  \frac{C}{\ell^3}\,\| f\|^2_{H^2(0,\ell)}.
\end{align*}
   \end{linenomath*} 

   Thus, the proof is complete. \hfill $\square$
   
\

   The main result in this section is the following:
   
\begin{teo}\label{uniqeta}
   Assume that $0 < \ell \leq L \leq L_*$, $0 < T_0 < T$, $u_x^\ell(0,t) = u_x^L(0,t)$ in~$(0,T)$, $\| u_0 \|_{L^2(0,L)} \leq M_0$ and~$|u_x^L(x,t)| \leq M$ in~$(0,\ell) \times (T_0,T)$. Then, there exists $\delta_0$ (only depending on~$L_*$, $T_0$, $T$, $M_0$ and~$M$) such that, if
   \begin{equation*}\label{cond-eta}
\int^{T}_{T_0} | \eta(t) |^2 \,dt \geq \delta_0,
   \end{equation*}
one necessarily has $\ell = L$.
\end{teo}

This result guarantees that uniqueness holds if the boundary input $\eta(t)$ is sufficiently large.

{\black
\noindent
\textbf{Proof:}
   The proof will be achieved by contradiction, assuming that $\ell < L$.
   
   Since $u^{\ell}(0,t) = u^{L}(0,t)$ and~$u_x^\ell(0,t) = u_x^L(0,t)$ in~$(T_0,T)$, we can use unique continuation to deduce that $u^{\ell}(x,t) = u^{L}(x,t)$ in~$(0,\ell) \times (T_0,T)$. Since $u^{\ell}(\ell,t) = 0$ in~$(0,T)$, we have that
   \begin{linenomath}
\begin{equation}\label{BC-uL}
u^{L}(\ell,t) = u^{L}(L,t) = 0 \ \text{ in } \ (0,T).
\end{equation}
   \end{linenomath}

   Let us introduce $\mathcal{D}(A) := H^2(\ell,L) \cap H^1_0(\ell,L)$ and~$Av = -v_{xx}$ for all $v \in \mathcal{D}(A)$ (the Dirichlet-Laplace operator in~$(\ell,L)$). Then,
   \[
\| u^L(\cdot\,,t)\|_{L^2(\ell,L)} = \| e^{-tA} u^L(\cdot\,,0) \|_{L^2(\ell,L)} \leq M_0 e^{-\zeta_1t} \ \ \forall t \in (T_0,T),
   \]
where $\zeta_1 > 0$ is the first eigenvalue of $A$, that is, $\zeta_1 = \pi^2(L-\ell)^{-2}$.

   In view of~\eqref{BC-uL}, we can deduce a standard energy estimate for $u^{L}_{xx}$ in~$(\ell,L) \times (0,T)$:
   \begin{linenomath*}
\begin{equation*}
\| u^{L}_{xx}(\cdot,t)\|_{L^2(\ell,L)}
= \| u^{L}_{t}(\cdot,t)\|_{L^2(\ell,L)}
\leq \frac{M_0}{t} e^{-\zeta_1t} \ \text{ in } \ (T_0,T).
\end{equation*}
   \end{linenomath*}
   Therefore,
\begin{linenomath*}
\begin{equation*}\label{expi}
\| u^{L}_{xx}(\cdot,t)\|_{L^2(\ell,L)}
\le \frac{M_0}{T_0}\exp\left( -\frac{\pi^2T_0}{(L-\ell)^2}\right) \ \text{ in } \ (T_0,T)
\end{equation*}
   \end{linenomath*}
and, in view of~Lemma~\ref{traceteo}, we have
   \begin{linenomath*}
\begin{equation}\label{expi2}
| u^{L}_x(\ell,t) | 
\le \frac{C(L_*) M_0}{T_0 (L-\ell)^{3/2}}
\exp\left( -\frac{\pi^2T_0}{(L-\ell)^2}\right) \ \text{ in } \ (T_0,T).
\end{equation}
   \end{linenomath*} 

   Taking the maximum of the right-hand side of \eqref{expi2} with respect to $L-\ell$, we see that
   \begin{linenomath*}
\begin{equation*}\label{frac1}
| u^{L}_x(\ell_1,t) | \le \frac{C(L_*,M_0)}{T_0^{7/4}} \ \text{ in } \ (T_0,T).
   \end{equation*}
   \end{linenomath*} 
   Since $u^{L}(\ell,t) = 0$ for $t>0$, we can apply Lemma~\ref{lemma2.8}  and get
   \begin{linenomath*}
\begin{equation}\label{tes}
\left( \int^{T}_{T_0} \int^{2\ep}_{\ep} | u^{L}(x,t) |^2 \,dx\,dt 
\right)^{\frac{1}{2}}
\le K(\ep,L_*,M_0,M) \frac{1}{T_0^{7\theta/4}} .
\end{equation}
   \end{linenomath*} 
   On the other hand, the Mean Value Theorem yields
   \begin{linenomath*}
\begin{equation*}
\int^{T}_{T_0} \int^{2\ep}_{\ep} \| u^{L}(x,t)\|^2 \,dx\,dt
= \ep \int^{T}_{T_0} \int^{2\ep}_{\ep} \| u^{L}(\xi_t,t)\|^2 \,dt ,
\end{equation*}
   \end{linenomath*} 
where $\xi_t \in (\ep, 2\ep)$ for all $t \in (T_0,T)$.
   Since $|u_x^L(x,t)| \leq M$ in~$(0,\ell) \times (T_0,T)$, we obtain:
   \begin{linenomath*}
\postdisplaypenalty=0
\begin{align*}
&\int^{T}_{T_0} \int^{2\ep}_{\ep} | u^{L}(x,t) |^2 \,dx\,dt
= \ep \int^{T}_{T_0} | u^{L}(0,t) + (u^{L}(\xi_t,t) - u^{L}(0,t)) |^2 \,dt
\\ & \quad \geq \frac{\ep}{2} \int^{T}_{T_0} | \eta(t) |^2 \,dt
- \frac{\ep}{2} \int^{T}_{T_0} | u^{L}(\xi_t,t) - u^{L}(0,t) |^2 \,dt
\\ & \quad \geq \frac{\delta_0 \ep}{2} - \frac{T M^2\ep^3}{2} .
\end{align*}
   \end{linenomath*} 
   Combining the last inequality and~\eqref{tes}, we see that
   \begin{linenomath*}
\begin{equation*}
\delta_0 \ep - T M^2 \ep^3 \leq 2K(\ep,L_*,M_0,M) \frac{1}{T_0^{7\theta/4}}
\end{equation*}
   \end{linenomath*} 
and, with $\ep$ small enough, the following holds:
   \begin{linenomath*}
\begin{equation}\label{contra}
\delta_0 \leq C(\ep,L_*,T,M_0,M) \frac{1}{\ep T_0^{7\theta/4}} .
\end{equation}
   \end{linenomath*} 

   Obviously, there exists (large) $T_0$ such that~\eqref{contra} is not satisfied.
   Thus, for this $T_0$, we reach a contradiction and we conclude that $\ell=L$.  \hfill $\square$

\begin{remark}
   When $\eta(t)\not\equiv 0$ but is small, in view of~Proposition~\ref{prop2.1}, we suspect that a non-uniqueness result can be obtained.
   Unfortunately,  to our knowledge this is unknown.
\end{remark}
}

\section{Results for the wave equation}\label{Sec-3}

   We will consider in this section the following problem for the 1D wave equation
   \begin{linenomath*}
\begin{equation}\label{problema-w}
\begin{cases}
u_{tt} - u_{xx} = 0,\, &0<x<\ell,\ \ 0 < t < T,\\
u(0,t)=\eta(t),\, u(\ell,t)=0,\, & 0 < t < T,\\
u(x,0)=u_0(x),\ \ \  u_t(x,0)=u_1(x),\, & 0<x<\ell,
\end{cases}
\end{equation}
   \end{linenomath*}
and the related inverse problem:

\

\textbf{Problem~IP-2 (wave equation, observation on the known boundary):}
   \textit{Assume that, in~\eqref{problema-w}, $u_0=u_0(x)$, $u_1=u_1(x)$ and $\eta=\eta (t)$ are given in appropriate spaces. Also, let the observation $u_x(0,\cdot)$ be known. Then, find $\ell$.}
   
\

   In this section, we assume that $u_0$ and $u_1$ are smooth and 
   \[
\frac{du_k}{dx}(0) = \frac{du_k}{dx}(\ell) = \frac{du_k}{dx}(L) = 0, \ \text{ for } \ k=0,1.
   \]
   Then, since $(u_k,\varphi_n) = O\left(\frac{1}{n^2}\right)$, the eigenfunction expansions of the solutions to~\eqref{problema-w} corresponding to~$\ell$ and~$L$ can be differentiated term-by-term with respect to~$x$ and the resulting series respectively converge in $C^0([0,\ell]\times [0,T])$ or $C^0([0,L]\times [0,T])$.

\

  We will analyze uniqueness:
   
\

\textbf{Uniqueness (observation on the known boundary):}
  \textit{Let $u^{\ell}$ and $u^L$ be the solutions to~\eqref{problema-w} respectively corresponding to~$\ell$ and~$L$ and assume that the corresponding observations coincide at~$x=0$, that is, $u^{\ell}_x(0,\cdot)=u^{L}_x(0,\cdot)$.
   Then, do we have $\ell=L$?}
   
\
   
   Since the equation is now hyperbolic and the information travels at finite speed, it will make sense to assume (at least) that $T > \ell$ in the sequel.

\subsection{Zero initial data and/or zero boundary data}\label{SSec-3.1}

   As before, some situations can be easily handled.

\subsubsection{ \textbf{Case I}: $\eta\not\equiv 0$ and $(u_0,u_1)\equiv(0,0)$}\label{SSSec-3.1.1}

   The argument is as in~Section~\ref{SSSec-2.1.1}:
   let us assume that~$0 < \ell < L < T$ and~$u^{\ell}_x(0,\cdot)=u^{L}_x(0,\cdot)$ and let us set $v := u^\ell - u^L$;
   then, by unique continuation, we have $v(x,t) = 0$ in~$(0,\ell) \times (\ell,T-\ell)$ and then, by energy estimates, $v(x,t) = 0$ in~$(0,\ell) \times (0,T)$, whence $u^L \equiv 0$ in~$(\ell,L) \times (0,T)$;
   but, using again unique continuation, this yields $\eta \equiv 0$, which is impossible.
   
   Consequently, we have:
   

\begin{prop}\label{prop3.1}
   Assume that $0 < \ell \leq L < T$, $\eta \in L^2(0,T)$ satisfies $\eta\not\equiv 0$ and $(u_0,u_1)\equiv(0,0)$.
   Then, if the solutions $u^\ell$ and $u^L$ respectively corresponding to~$\ell$ and~$L$ satisfy $u^{\ell}_x(0,\cdot)=u^{L}_x(0,\cdot)$, we necessarily have $\ell = L$.
\end{prop}

\subsubsection{ \textbf{Case II}: $\eta\equiv 0$ and $(u_0,u_1)\not\equiv(0,0)$}\label{SSSec-3.1.2} 

   As in Section~\ref{SSSec-2.1.2}, we can easily obtain here explicit expressions of the solutions corresponding to $\ell$ and~$L$.
   Thus, with the notation used there, one has:
\begin{linenomath*}
\begin{align*}
u^{\ell}(x,t)=& \sum_{n=1}^{\infty} 
\left[ (u_0, \va_n)_{\ell} \cos(\sqrt{\lambda_n} t)
+ \dfrac{1}{\sqrt{\lambda_n}} (u_1, \va_n)_{\ell} \sin(\sqrt{\lambda_n} t) \right]
\va_n(x) ,
\end{align*}
\end{linenomath*}
\begin{linenomath*}
\begin{align*}
u^{L}(x,t)=& \sum_{n=1}^{\infty} 
\left[ (u_0, \psi_n)_{L} \cos(\sqrt{\mu_n} t)
+ \dfrac{1}{\sqrt{\mu_n}} (u_1, \psi_n)_{L} \sin(\sqrt{\mu_n} t) \right]
\psi_n(x) .
\end{align*}
\end{linenomath*}

   It is easy to find arbitrary close 
   (but different) $\ell$ and~$L$ and couples $(u_0,u_1)$
   (defined for all~$x \in \R_+$) such that the corresponding observations coincide.
   
   More precisely, the following holds:
   

\begin{prop}\label{prop3.2}
   Assume that $\eta \equiv 0$, $0 < \ell < L$ and $L/\ell = m_0/n_0 \in \Q$ and let us set
   \begin{equation*}\label{eq3.2a-w}
u_0(x) := \sin \dfrac{n_0 \pi x}{\ell} \equiv \sin \dfrac{m_0 \pi x}{L}, \quad u_1 \equiv 0.
   \end{equation*}
   Then, the corresponding $u^\ell$ and $u^L$ satisfy~$u^{\ell}_x(0,\cdot)=u^{L}_x(0,\cdot)$.
   The same result holds for
    \begin{equation*}\label{eq3.2b-w}
u_0(x) \equiv 0, \quad u_1(x) := \sin \dfrac{n_0 \pi x}{\ell} \equiv \sin \dfrac{m_0 \pi x}{L}.
   \end{equation*}
\end{prop}

\begin{remark}
   Since we can find $\ell$ and $L$ with rational $L/\ell$ and arbitrarily small $L-\ell$, this proposition indicates that non-uniqueness holds even if $L$ is assumed to be close to $\ell$.
    On the other hand, if we choose a nonzero boundary input $\eta$, we can deduce an upper estimate of~$\vert L-\ell\vert$ where the bound becomes smaller as $\eta$ becomes larger.
   This is proved in the next section.
 \end{remark}

\subsection{An asymptotic result for nonzero initial and boundary data}\label{SSec-3.2}

{\black
   Let us assume that $\eta \not\equiv 0$ and $(u_0,u_1) \not\equiv (0,0)$. 
If we try to prove uniqueness for sufficiently large $\eta$ by arguing as for 
the heat equation, then a difficulty arises. Indeed, we would need decay in time of the solution; but this does not hold for the classical wave equation. Note however that something can be said for instance for the {\it telegraphist-Klein-Gordon} equation $u_{tt} + a(x)u_t - u_{xx} + m(x) u = 0$ with $a, m \geq 0$.

   In this section, we will prove an asymptotic uniqueness result. To this purpose, we will need the following technical result:

\begin{lem}\label{lemma1-new}
   Let $w = w(x,t)$ satisfy
   \begin{equation}\label{23p}
\left\{ \begin{array}{ll}
w_{tt} - w_{xx} = 0, \quad &\ell<x<L, \ \  T_0 < t < T_1, \\
w(\ell,t) = w(L,t) = 0, \quad & T_0 < t<T_1. \\
\end{array}\right.
   \end{equation}
Then
   \begin{equation}\label{23pp}
\| w_x(\ell,\cdot)\|_{L^2(T_0,T_1)}^2
\le \frac{T_1-T_0 + 2(L-\ell)}{L-\ell} 
\left( \| w_t(\cdot\,,T_0) \|_{L^2(\ell,L)}^2 + \| w_x(\cdot\,,T_0) \|_{L^2(\ell,L)}^2 \right).
   \end{equation}
\end{lem}

\noindent
{\bf Proof:} This is called a direct inequality for $w$: it allows to estimate the ``lateral'' flux $w_x(\cdot\,,T_0)$ in terms of the energy associated to~$w$ at time~$T_0$.

   It can be proved by the {\it multiplier method.} Indeed, the change of variables $\xi = x-\ell$, $\tau = t-T_0$ leads to the following reformulation of~\eqref{23p}:
   \[
\left\{ \begin{array}{ll}
z_{\tau\tau} - z_{\xi\xi} = 0, \quad &0 < \xi < d, \ \  0 < \tau < T_*, \\
z(0,\tau) = z(d,\tau) = 0, \quad & 0 < \tau <T_*, \\
\end{array}\right.
   \]
where $d = L - \ell$ and~$T_* = T_1 - T_0$. For any $q \in C^1([0,d])$, we can directly prove the identity
   \[
\frac{1}{2} \int^{T_*}_0 \int^d_0 q_\xi (|z_\tau|^2 + |z_\xi|^2) \,d\xi\,d\tau 
+ \Bigl. \int^d_0 q(\xi)  z_\xi z_\tau \,d\xi \Bigr|^{\tau=T_*}_{\tau=0}
= \Bigr. \frac{1}{2} \int^{T_*}_0 q(\xi) |z_\xi|^2 \,d\tau \Bigr|^{\xi=d}_{\xi=0} .
   \]
   Taking $q(\xi) \equiv \xi-d$ and using energy conservation, we get:
   \[
\begin{array}{l} \dis
\frac{d}{2} \int^{T_*}_0 q(\xi) |z_\xi(0,\tau)|^2 \,d\tau
= \frac{T_*}{2} \int^d_0 (|z_\tau(\xi,0)|^2 + |z_\xi(\xi,0)|^2) \,d\xi 
+ \Bigl. \int^d_0 (\xi-d)  z_\xi z_\tau \,d\xi \Bigr|^{\tau=T_*}_{\tau=0}
\\ \dis \phantom{\frac{d}{2} \int^{T_*}_0 q(\xi) |z_\xi(0,\tau)|^2 \,d\tau}
\leq \frac{1}{2}(T_* + 2) \int^d_0 (|z_\tau(\xi,0)|^2 + |z_\xi(\xi,0)|^2) \,d\xi ,
\end{array}
   \]
whence we readily obtain~\eqref{23pp}.  \hfill $\square$

\

   In the sequel, we will denote by $E_\ell = E_\ell(t)$ the energy associated to the solution to~\eqref{problema-w} in~$(0,\ell)$, that is,
   \[
E_\ell(t) := \int^{\ell}_0 \left( | u^\ell_t(x,t) |^2 + | u^\ell_x(x,t) |^2 \right) \,dx.
   \]

The main result in this section is estimation  of lengths of spatial
intervals which can be more accurate when amplitudes of the boundary 
inputs are larger:

\begin{teo}\label{asymp-uniq}
   Assume that $0 < \ell_0 \leq \ell \leq L \leq \ell_1$, $T > 4 \ell_1$, $u_x^\ell(0,t) = u_x^L(0,t)$ in~$(0,T)$ and~$E_\ell(t), E_L(t) \leq M$ in~$(0,T)$. Then, there exists $\delta_0>0$ such that
   \begin{equation}\label{2.2p}
\sup_{t \in [2\ell, T-2\ell]} | \eta(t) | \geq \delta_0 \ \Rightarrow \ 
| L-\ell| \leq \frac{\ell_1}{2}\frac{(T+2\ell_1-4\ell_0)M}{\delta_0^2} \,.
   \end{equation}
\end{teo}

\noindent
{\bf Proof:}
   Let us set $w = u^\ell - u^L$. Then, from unique continuation, we have
   \[
w(\ell,t) = 0 \ \text{ in } \ (\ell,T-\ell),
   \]
whence $u^L$ satisfies:
   \[
\left\{ \begin{array}{l}
u^L_{tt} - u^L_{xx} = 0, \quad \ell<x<L, \ \  \ell<t<T-\ell,
\\ \noalign{\vskip 3pt}
u^L(\ell,t) = u^L(L,t) = 0, \quad \ell<t<T-\ell.
\end{array}\right.
   \]
   Hence, in view of~Lemma~\ref{lemma1-new}, the following holds:
   \begin{equation}\label{3.2}
\|u^L_x(\ell,\cdot)\|_{L^2(\ell,T-\ell)}
\leq \sqrt{ \frac{(T-2\ell) + 2(L-\ell)M}{L-\ell} }
\leq \frac{\sqrt{ (T+2\ell_1-4\ell_0)M }}{\sqrt{L-\ell}}.
  \end{equation}
   Since $2\ell_0 < t_0 < T-2\ell_0$ and~$\ell, L \in [\ell_0,\ell_1]$, one has $2\ell < T-2\ell$ and for any~$t \in [2\ell,T-2\ell]$, D'Alembert formula yields
   \[
u^L(0,t) = -\frac{1}{2}\int^{t+\ell}_{t-\ell} u^L_x(\ell,s) \,ds.
   \]
   Consequently, from the Cauchy-Schwarz inequality and~\eqref{3.2}, we find that
   \[
| u^L(0,t) |^2 \leq \frac{\ell}{2} \| u^L_x(\ell,\cdot) \|_{L^2(\ell,L-\ell)}^2 
\leq \frac{\ell_1}{2} \frac{(T+2\ell_1-4\ell_0)M}{L-\ell}
   \]
and, if $\displaystyle{\sup_{t \in [2\ell, T-2\ell]} | \eta(t) | \geq \delta_0}$, one has
   \[
\delta_0^2 \leq \frac{\ell_1}{2} \frac{(T+2\ell_1-4\ell_0)M}{L-\ell}\,.
   \]
   This ends the proof. \hfill $\square$
   
\

\begin{remark}
   Obviously, the counter-example in Proposition~\ref{prop3.2} shows that, in general, uniqueness cannot be expected.
   However, there are several situations where the estimate in~\eqref{2.2p} can be useful:
   if~$\delta_0$ is large with respect to the energy bound~$M$ or~$\ell_1 \sim 0$, or $T\sim 2\ell_0$ and $\ell_0\sim \ell_1$, etc.
\end{remark}

\begin{remark}
If $u_1$ vanishes, then we can improve the previous result. Indeed, 
the solutions $u^\ell$ and $u^L$ can be extended as even functions to negative 
times and we can take $M = \| u_{0,x} \|_{L^2}^2$. For brevity, we omit the 
details.
\end{remark}
}

\section{Some numerical results}\label{Sec-4}

   In this section, we will consider the previous inverse problems for the heat and wave equations. As usual, we will carry out the reconstruction of the unknown length through the resolution of some appropriate extremal problems.
   This strategy has been applied in some previous papers of the authors for other similar problems, see~\cite{DouFe}--\cite{DouFe3}.
   The results of the numerical tests that follow will serve to illustrate the theoretical results in the previous sections. 

\subsection{Tests 1: The classical heat equation}\label{SSec-4.1}

   We deal with the following
   
\
   
   \textbf{Reformulation of IP-1:} \textit{Given $\eta = \eta(t)$, $u_0 = u_0(x)$, $T>0$ and $\beta = u_x(0,\cdot)$, find $\ell \in (\ell_0,\ell_1)$ such that 
   \begin{equation*}\label{pb.opt}
J(\ell) \le  J(\ell'),\quad \forall\, \ell'\in (\ell_0,\ell_1) ,
   \end{equation*}
where $J$ is given by
   \begin{equation}\label{eq.J}
J(\ell) = \dfrac{1}{2}\displaystyle\int_0^T |\beta(t) - u^\ell_x (0,t)|^2\, dt.
   \end{equation}
   Here, $u^\ell$ is the {\it state,} i.e.\ the solution to~\eqref{problema}, corresponding to the unknown length~$\ell$.}
   
\
   
   Three different situations will be analyzed for the heat equation.
   In the first two cases, we will check that uniqueness holds: zero initial data and nonzero initial data and sufficiently large $\eta$. In the third case we will consider a non-uniqueness situation corresponding to some nonzero initial data and ``small'' $\eta$ and we will study the behavior of the numerical  algorithm.
   To this purpose, we will implement the \texttt{fmincon} function from the MatLab Optimization Toolbox using the \texttt{active-set} minimization algorithm. 

\

\noindent
\textbf{Case 1.1: Heat equation with $u_0=0$ and $\eta\neq 0$.}

   We take $T=5$,  $\eta(t) = 5\sin^3 t$ in~$(0,T)$ and~$u_0(x) \equiv 0$.
   Starting from $L_{i} = 3$, our goal is to recover the desired value of the length $L_d= 2$.

\begin{figure}[h!]
\begin{minipage}[t]{0.49\linewidth}
\noindent
\includegraphics[width=\linewidth]{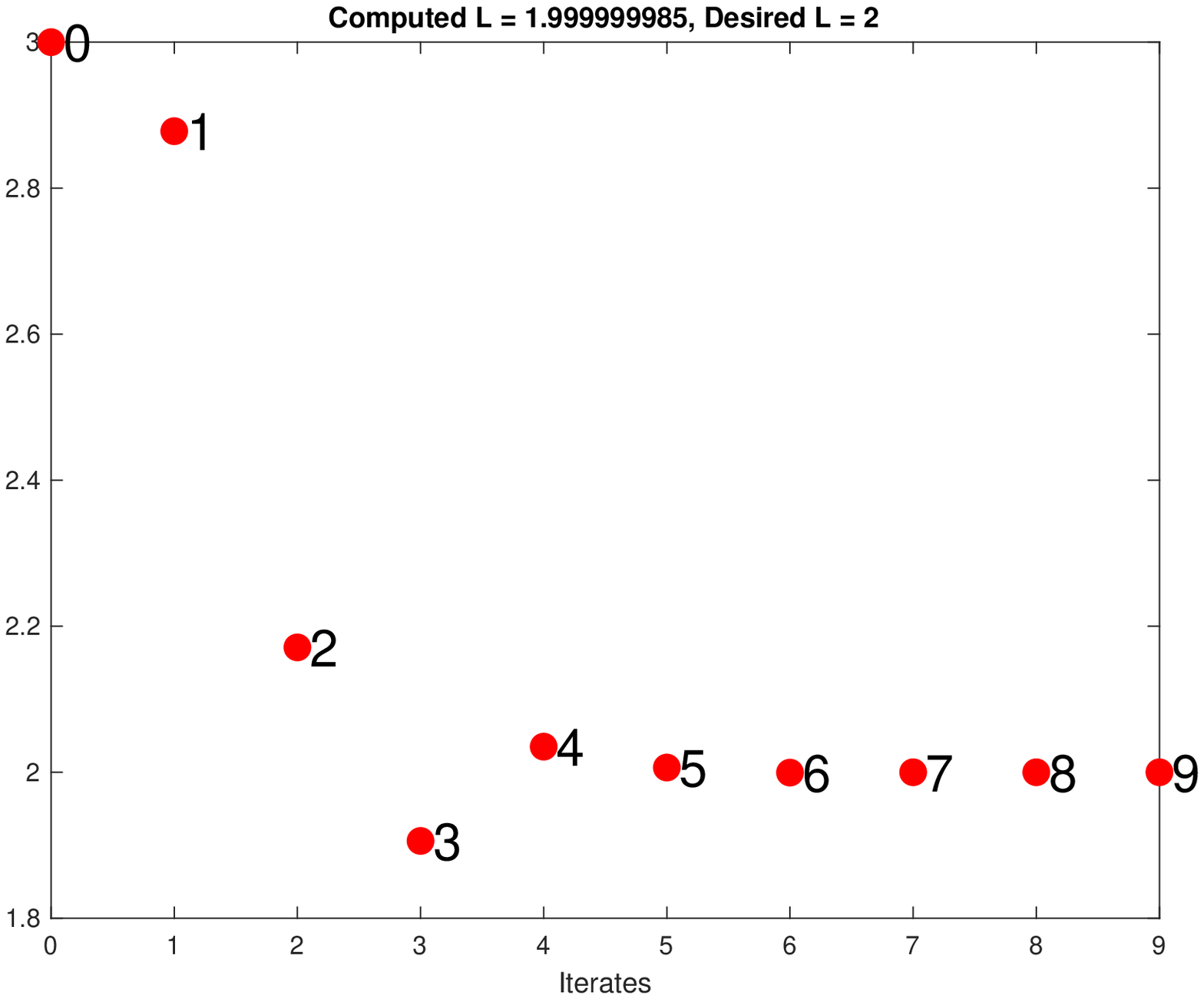}
\caption{Heat equation with $u_0=0$ and $\eta\neq 0$. The iterates in \texttt{active-set} algorithm.}
\label{Heat_iter_case1}
\end{minipage}
\hfill
\begin{minipage}[t]{0.49\linewidth}
\noindent
\includegraphics[width=\linewidth]{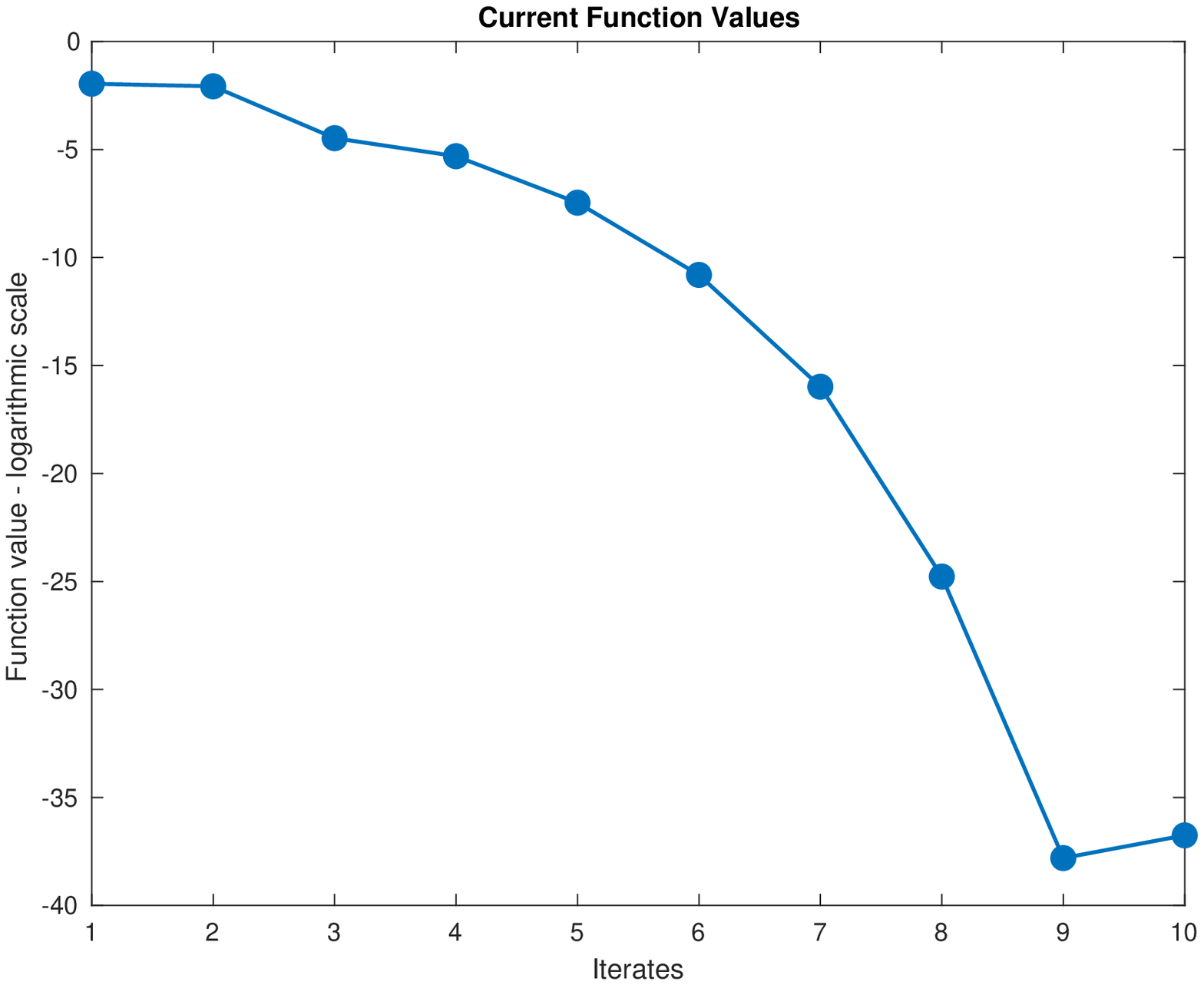}
\caption{Heat equation with $u_0=0$ and $\eta\neq 0$. Evolution of the cost.}
\label{Heat_J_case1}
\end{minipage}
\end{figure}

   The results of this numerical experiments can be seen in~Table~\ref{tab:heat_case1}, where the effect of random noise in the target are shown.
   The computed length is denoted by~$L_c$.
   The corresponding solution to~\eqref{problema} is displayed in~Figure~\ref{Heat_Sol_case1}.
   The evolution of the iterates and the cost in the minimization process in the absence of the random noise appear in~Figures~\ref{Heat_iter_case1} and~\ref{Heat_J_case1}, respectively.
   

%
\begin{table}[h!]
\begin{minipage}[t]{0.46\linewidth}
\centering
\renewcommand{\arraystretch}{1.2}
\medskip
\caption{Heat equation - $u_0=0$ and $\eta\neq 0$. Results with random noise in the target
{(the desired length is $L_d = 2$). } \\
{ }\\
{ }}
\begin{tabular}{|cccc|} \hline
\%  noise  & Cost  & Iterates &  Computed $L_c$ \\
\hline
1\%           &   1.e-4         &       10       &   1.997586488 \\ \hline
0.1\%        &    1.e-6        &      9        &   1.999864829\\ \hline
0.01\%      &   1.e-9       &       8         &   2.000017283\\\hline
0.001\%    &   1.e-1       &     9         &   1.999998535\\ \hline
0\%           &   1.e-16       &     9         &    1.999999991  \\  \hline
\end{tabular}
\label{tab:heat_case1}
\end{minipage}
\hfill
\begin{minipage}[t]{0.54\linewidth}
\vspace{0cm}
\includegraphics[width=\linewidth]{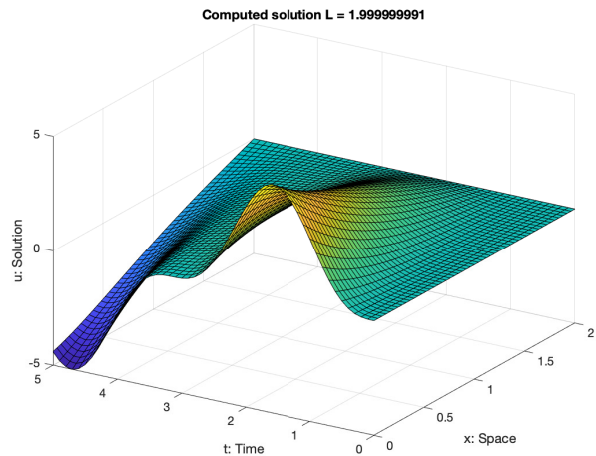}
\captionof{figure}{Heat equation with $u_0=0$ and $\eta\neq 0$. The computed solution.}
\label{Heat_Sol_case1}
\end{minipage}
\end{table}
%

\eject

\

\noindent 
\textbf{Case 1.2: Heat equation with $u_0\neq 0$ and large $\eta$.}

   We take $T=5$,  $\eta(t) = 0.2 \, t \, (2+t)$ in~$(0,T)$ and~$u_0(x) \equiv 5 x (2-x)$.
   Now, starting from $L_{i} = 0.5$, the target value that we want to recover is $L_d= 2$.
   
\begin{figure}[h!]
\begin{minipage}[t]{0.49\linewidth}
\noindent
\includegraphics[width=\linewidth]{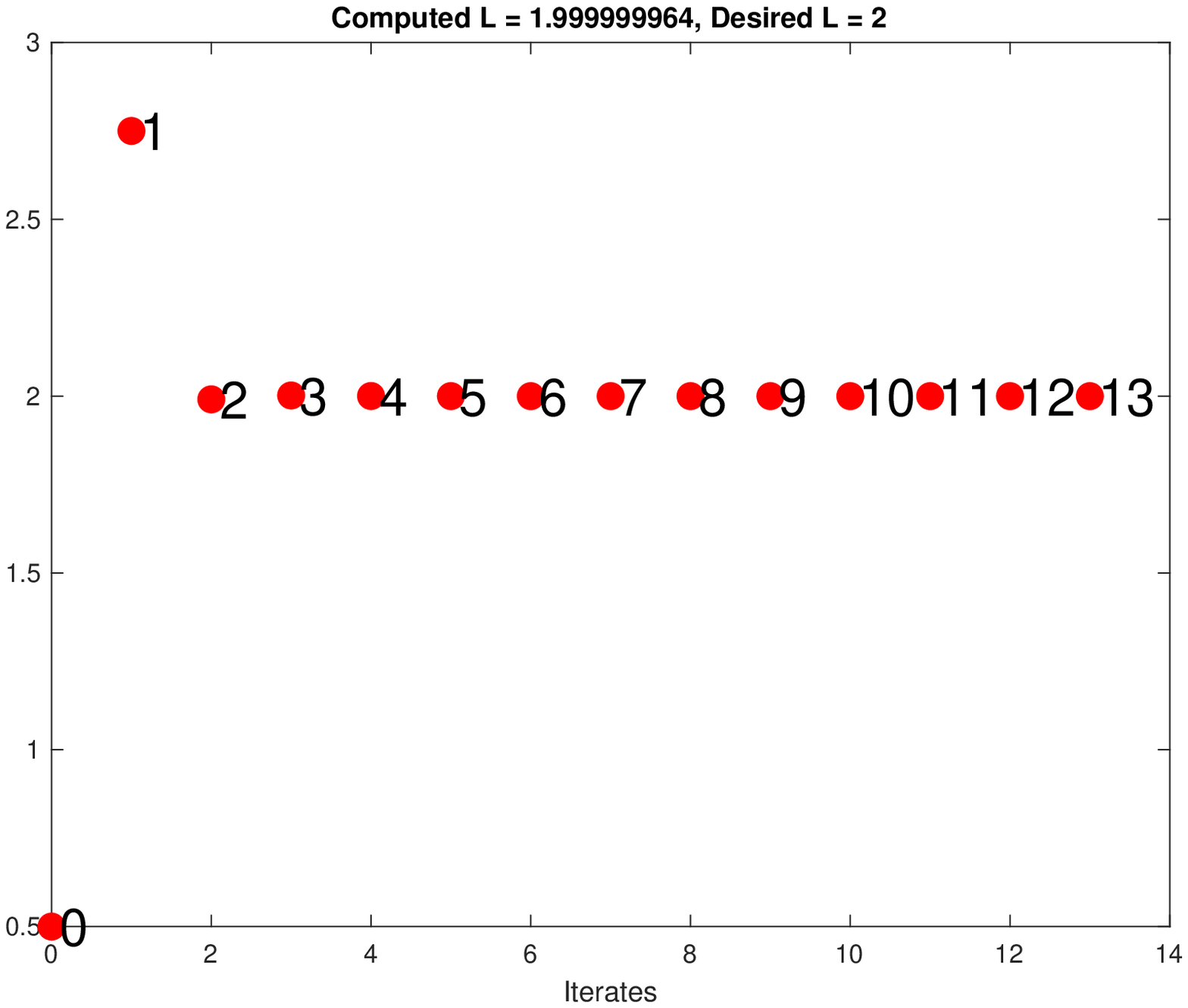}
\caption{Heat equation, fixed~$u_0$ and~large~$\eta$. The iterates in~\texttt{active-set} algorithm.}
\label{Heat_iter_case2}
\end{minipage}
\hfill
\begin{minipage}[t]{0.49\linewidth}
\noindent
\includegraphics[width=\linewidth]{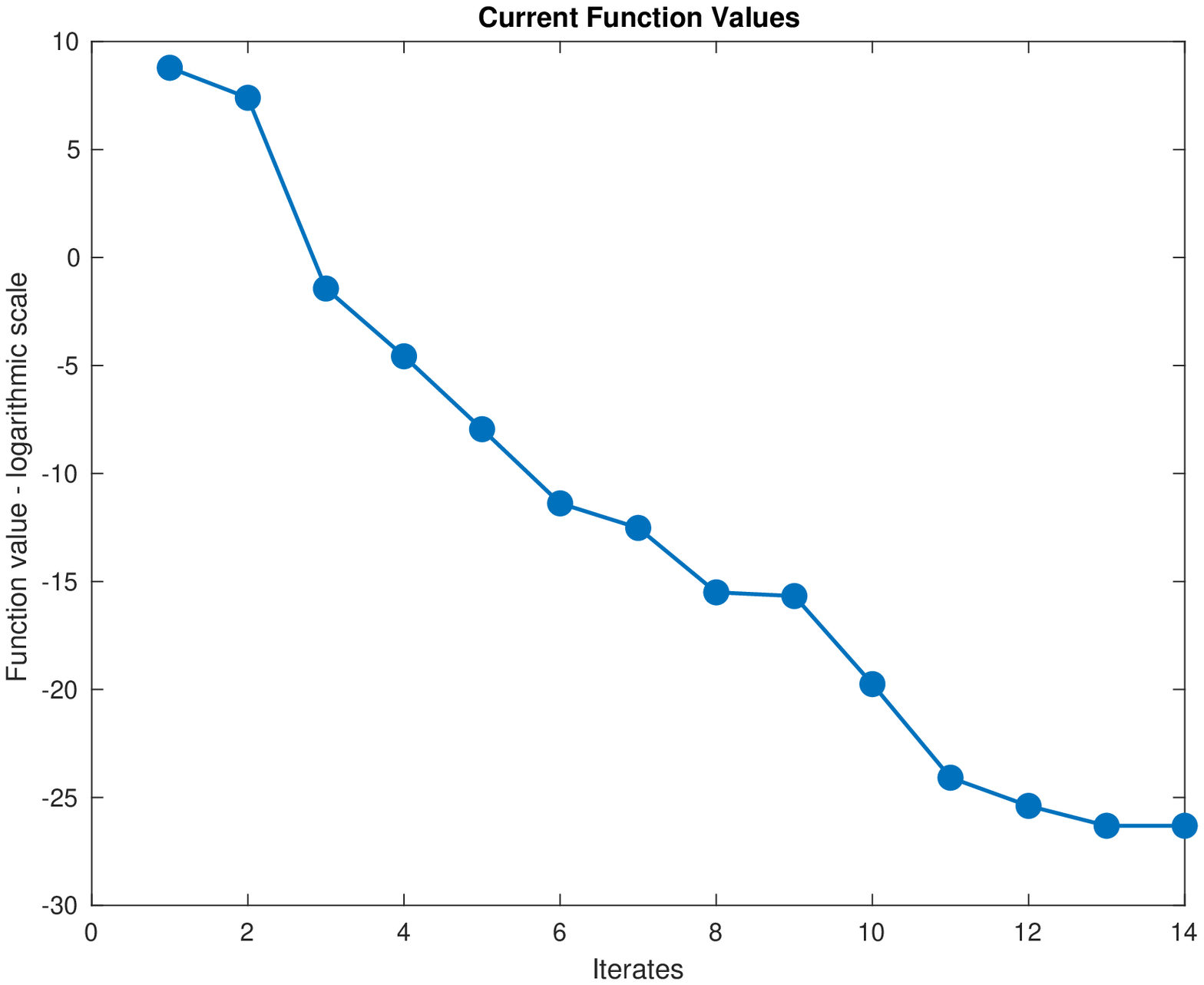}
\caption{Heat equation, fixed~$u_0$ and~large~$\eta$. Evolution of the cost.}
\label{Heat_J_case2}
\end{minipage}
\end{figure}

   The results of the numerical implementation are shown in Table~\ref{tab:heat_case2}, where again random noise was incorporated. The contents of Figures~\ref{Heat_Sol_case2}, \ref{Heat_iter_case2} and~\ref{Heat_J_case2} are similar to those above.
      

%
\begin{table}[h!]
\begin{minipage}[t]{0.46\linewidth}
\centering
\renewcommand{\arraystretch}{1.2}
\medskip
\caption{Heat equation - fixed $u_0$ and large~$\eta$. Results with random noise in the target
{(the desired length is $L_d = 2$). } \\
{ }
 \\
{ }}
\begin{tabular}{|cccc|} \hline
\%  noise  & Cost  & Iterates &  Computed $L_c$ \\
\hline
1\%           &   1.e-3         &       8       &   1.999952948  \\ \hline
0.1\%        &    1.e-5        &      11        &  2.000008948\\ \hline
0.01\%      &   1.e-8       &       9          &   2.000003174\\\hline
0.001\%    &   1.e-10       &     7           &   1.999999932\\ \hline
0\%           &   1.e-12       &     13         &    1.999999964  \\  \hline
\end{tabular}
\label{tab:heat_case2}
\end{minipage}
\hfill
\begin{minipage}[t]{0.54\linewidth}
\vspace{0cm}
\includegraphics[width=\linewidth]{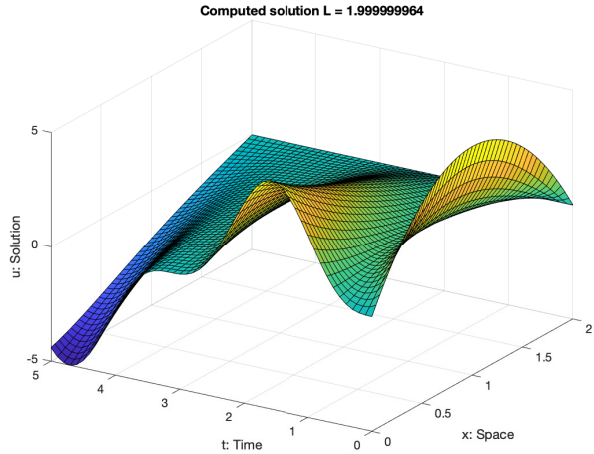}
\captionof{figure}{Heat equation, $u_0\neq 0$ and~large~$\eta$. The computed solution.}
\label{Heat_Sol_case2}
\end{minipage}
\end{table}
%

\

\eject
{\black
\noindent 
\textbf{Case 1.3: Heat equation with $u_0\neq 0$ and ``small'' $\eta$.}}

   Here, we deal with a non-uniqueness situation.
   Our aim is to investigate the behavior of the algorithm in a situation of this kind.
   
   We take $T=5$, $\eta=0$ in~$(0,T)$ and~$u_0(x) \equiv \sin(\pi x/2)$.
   Note that we have $u_0(x) \equiv \sin(3\pi x/L^1_d) \equiv \sin(2\pi x/L^1_d)$, with $L^1_d = 6$ and~$L^2_d = 4$; consequently, this initial data can be used as in~Section~\ref{SSSec-2.1.2} to prove non-uniqueness.
   
\begin{figure}[h!]
\begin{minipage}[t]{0.47\linewidth}
\noindent
\includegraphics[width=\linewidth]{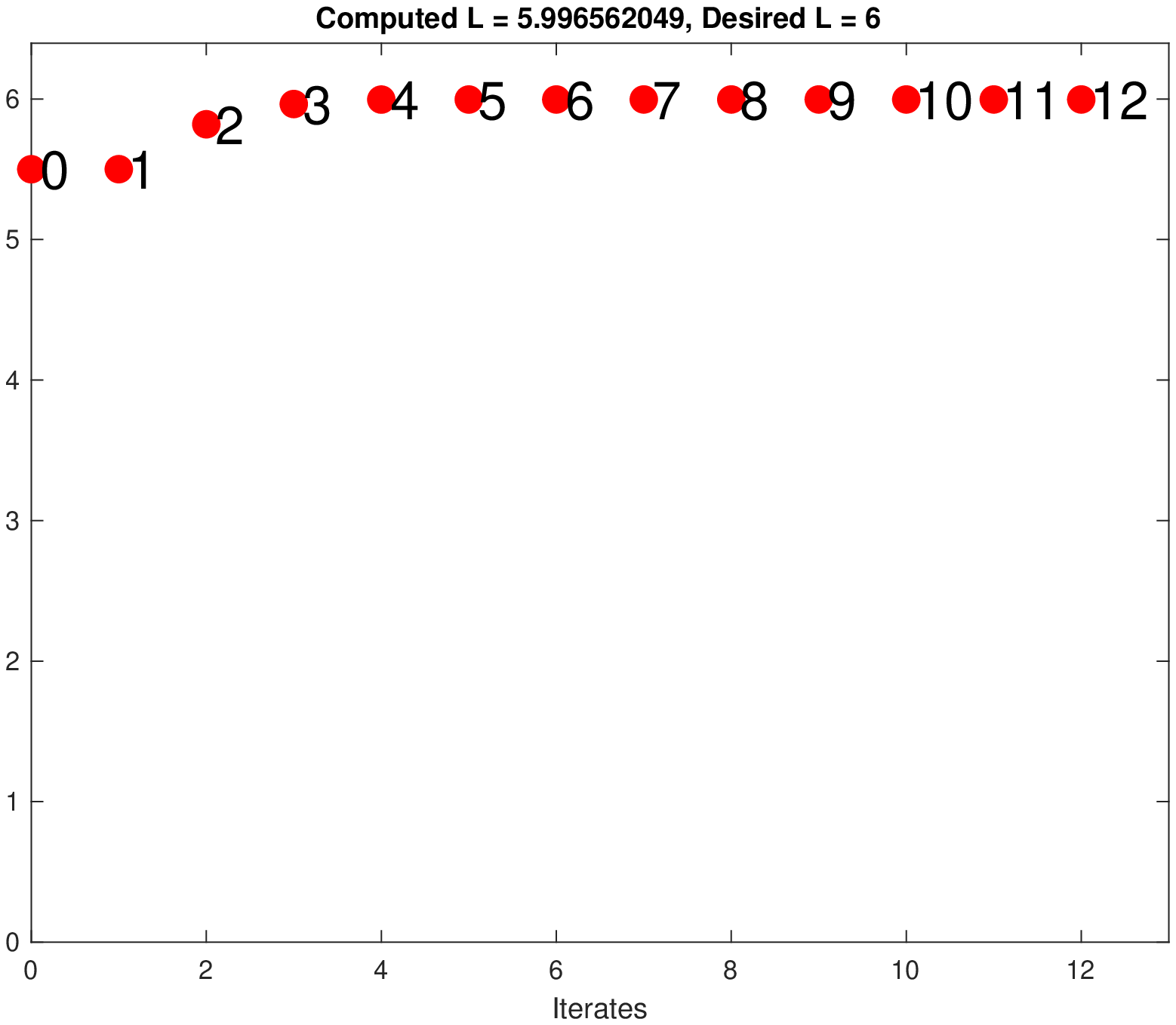}
\caption{Heat equation - $\eta=0$,  fixed $u_0(x)$. Iterates in \texttt{active-set} algorithm with $L^1_d=6$.}
\label{Heat_iter_case3}
\end{minipage}
\hfill
\begin{minipage}[t]{0.49\linewidth}
\noindent
\includegraphics[width=\linewidth]{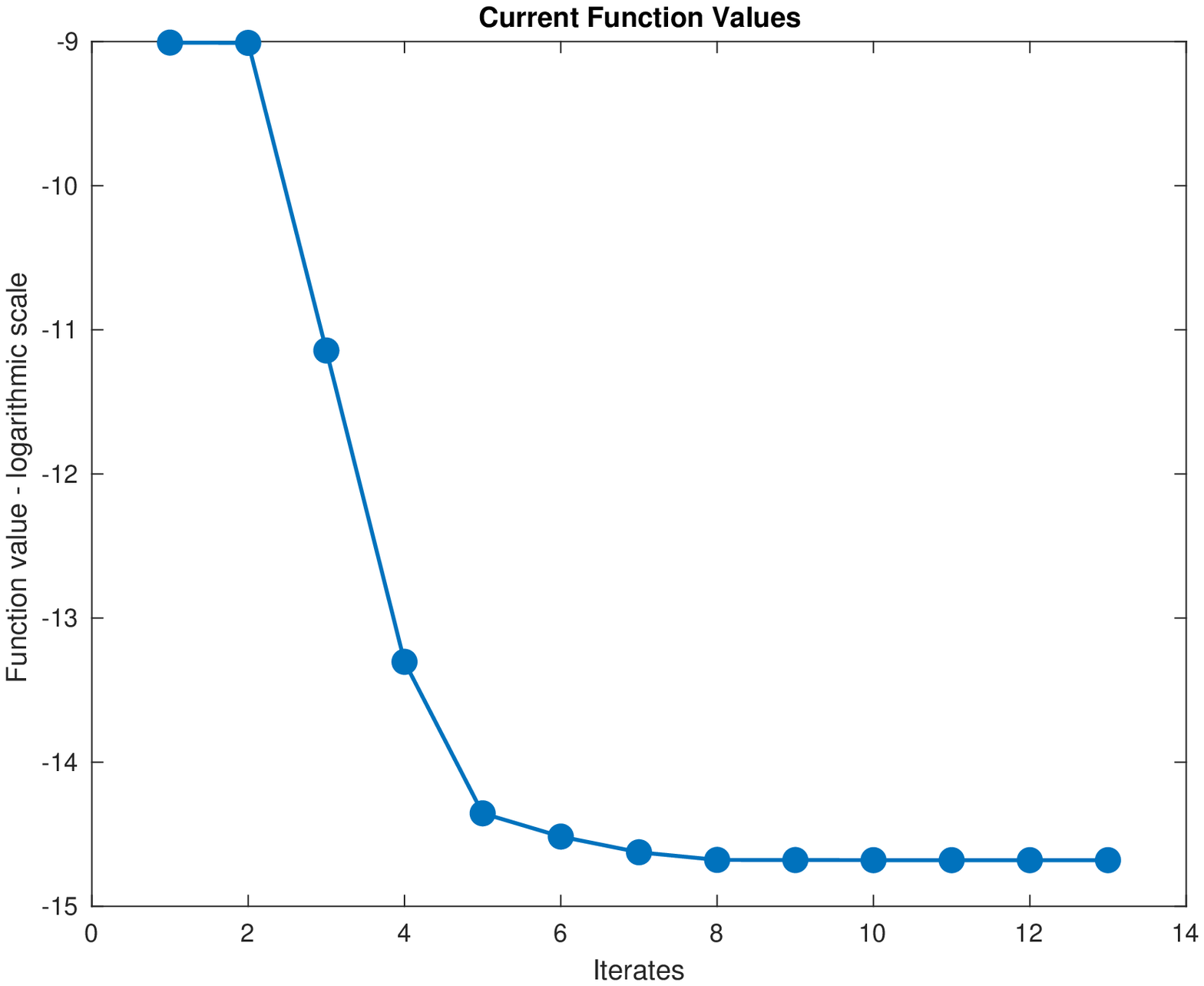}
\caption{Heat equation - $\eta=0$,  fixed $u_0(x)$. Evolution of the cost for $L^1_d=6$, $J(L^1_c)<10^{-7}$.}
\label{Heat_J_case3}
\end{minipage}
\end{figure}

\begin{figure}[h!]
\begin{minipage}[t]{0.49\linewidth}
\noindent
\includegraphics[width=\linewidth]{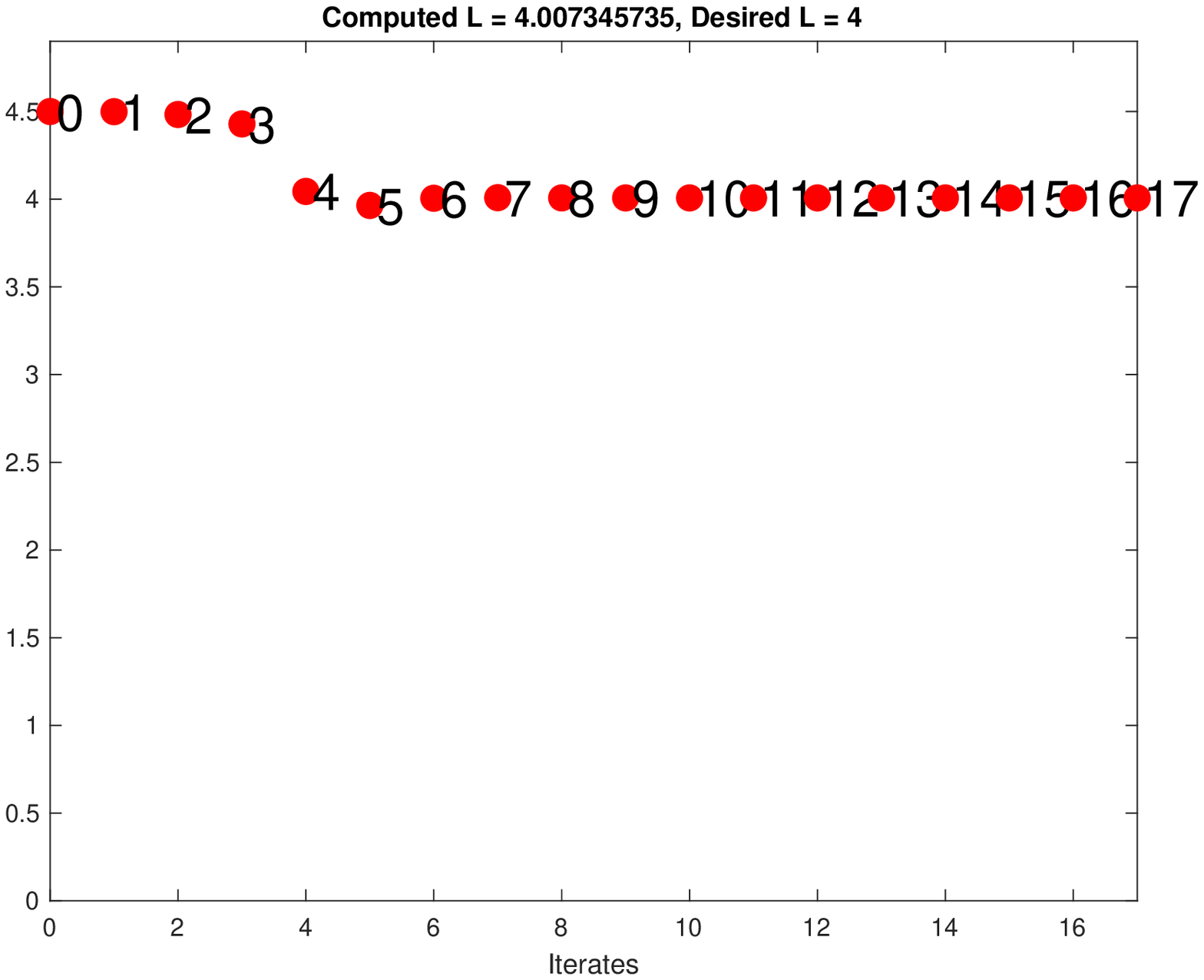}
\caption{Heat equation -  $\eta=0$,  fixed $u_0(x)$. Iterates in \texttt{active-set} algorithm with $L^2_d = 4$.}
\label{Heat_iter_case3b}
\end{minipage}
\hfill
\begin{minipage}[t]{0.49\linewidth}
\noindent
\includegraphics[width=\linewidth]{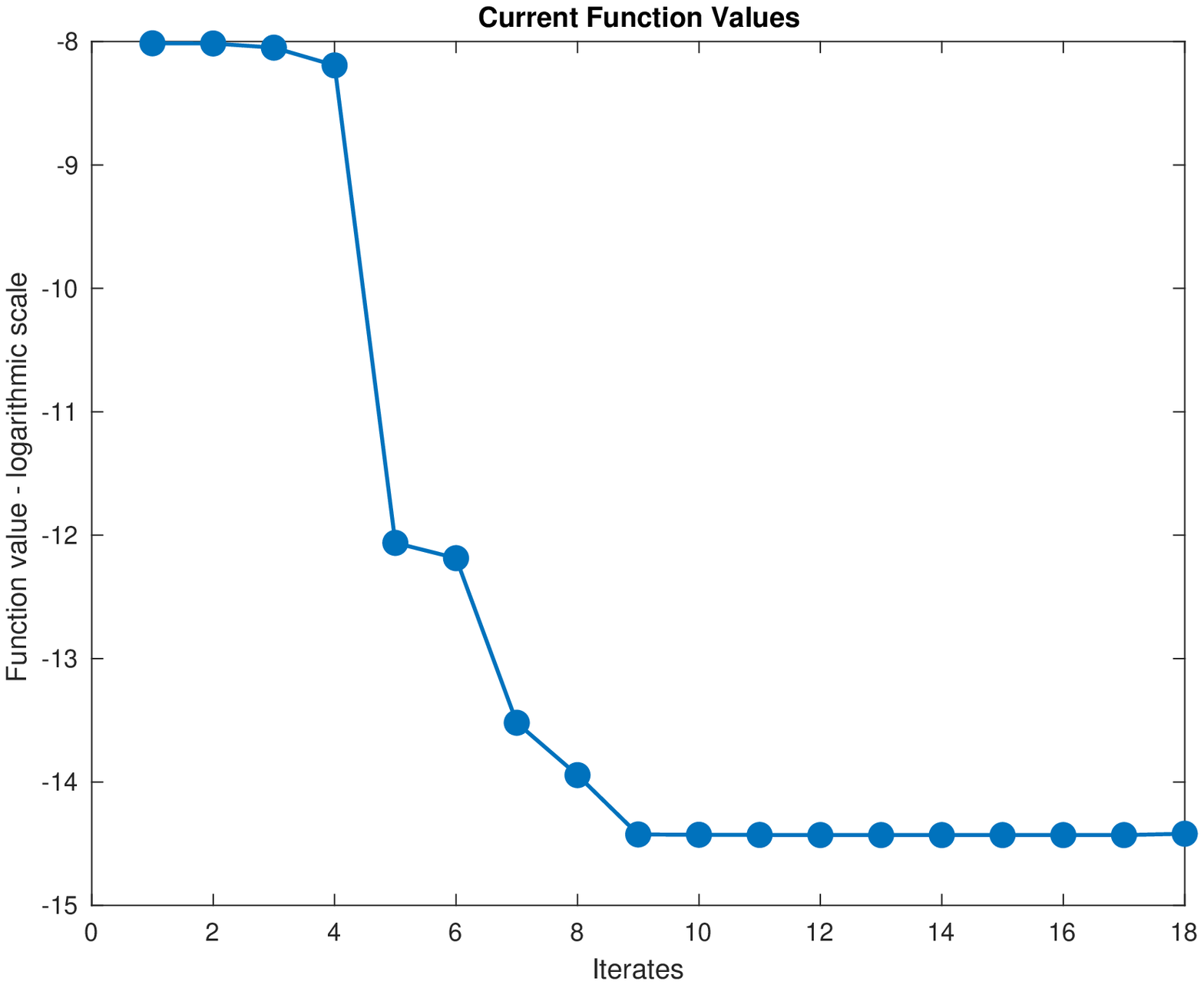}
\caption{Heat equation - $\eta=0$,  fixed $u_0(x)$. Evolution of the cost for $L^2_d=4$, $J(L_c^2)<10^{-7}$.}
\label{Heat_J_case3b}
\end{minipage}
\end{figure}

   We will consider the following experiments:
   
\begin{itemize}

\item First, we start {\black from~$L_i = 5.5$,} and we obtain the results exhibited in~Figures~\ref{Heat_iter_case3} and~\ref{Heat_J_case3}.
   The computed value is $L^1_c=5.996562049$ and the associated cost is $J(L_c^1)< 10^{-7}$.
   
\item Then, we start {\black from~$L_i = 4.5$,} and we obtain the results exhibited in~Figures~\ref{Heat_iter_case3b} and~\ref{Heat_J_case3b}.
   The computed value is $L^1_c=4.007345905$ and the associated cost is again $J(L_c^1)< 10^{-7}$.
   
\end{itemize}

\begin{figure}[h!]
\begin{minipage}[t]{0.49\linewidth}
\noindent
\includegraphics[width=\linewidth]{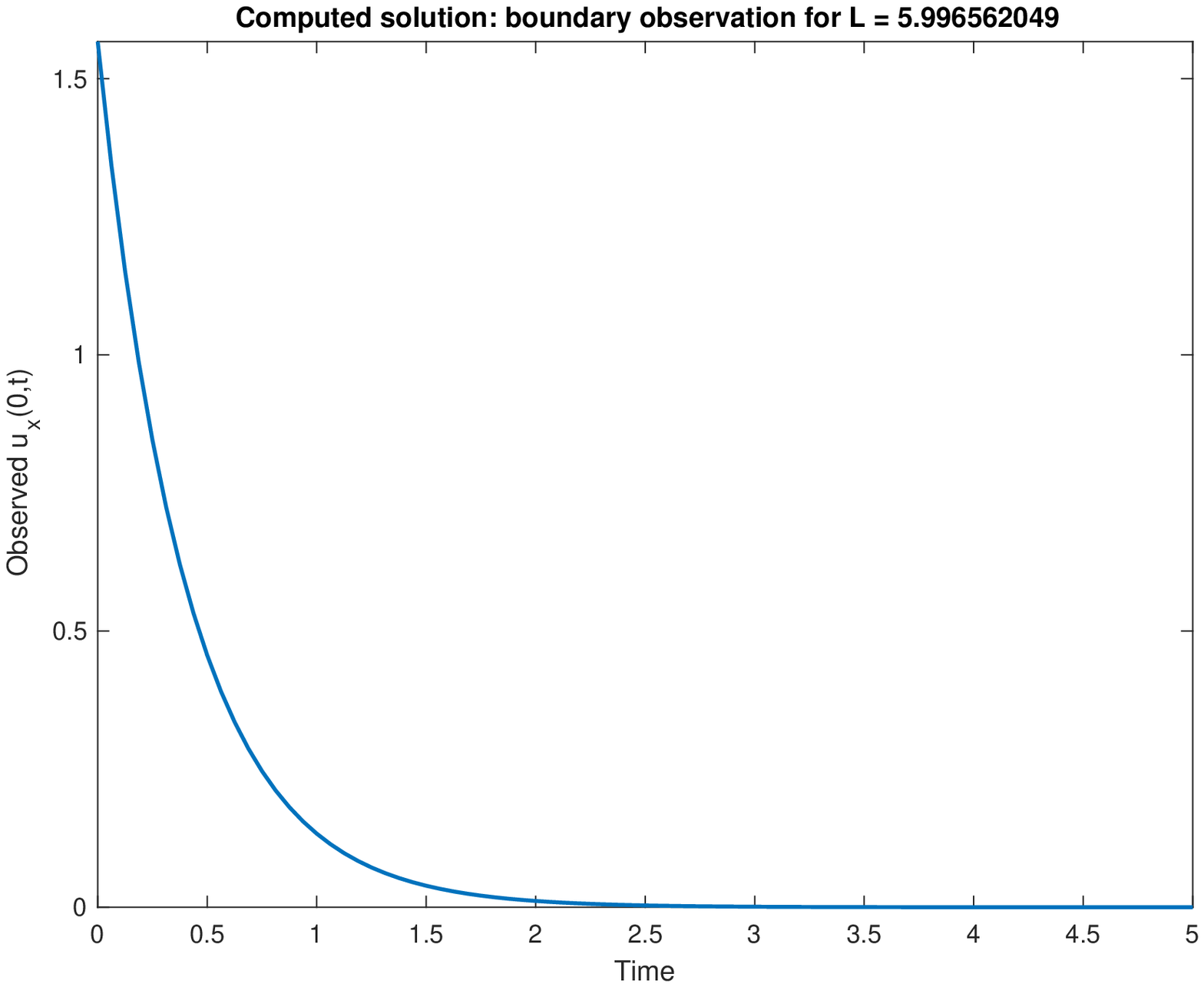}
\caption{Heat equation -  fixed $u_0$ and~$\eta=0$. The computed boundary observation $u_x(0,\cdot)$ for~$L^1_c=5.996562049$.}
\label{Heat_Sol_case3}
\end{minipage}
\hfill
\begin{minipage}[t]{0.49\linewidth}
\noindent
\includegraphics[width=\linewidth]{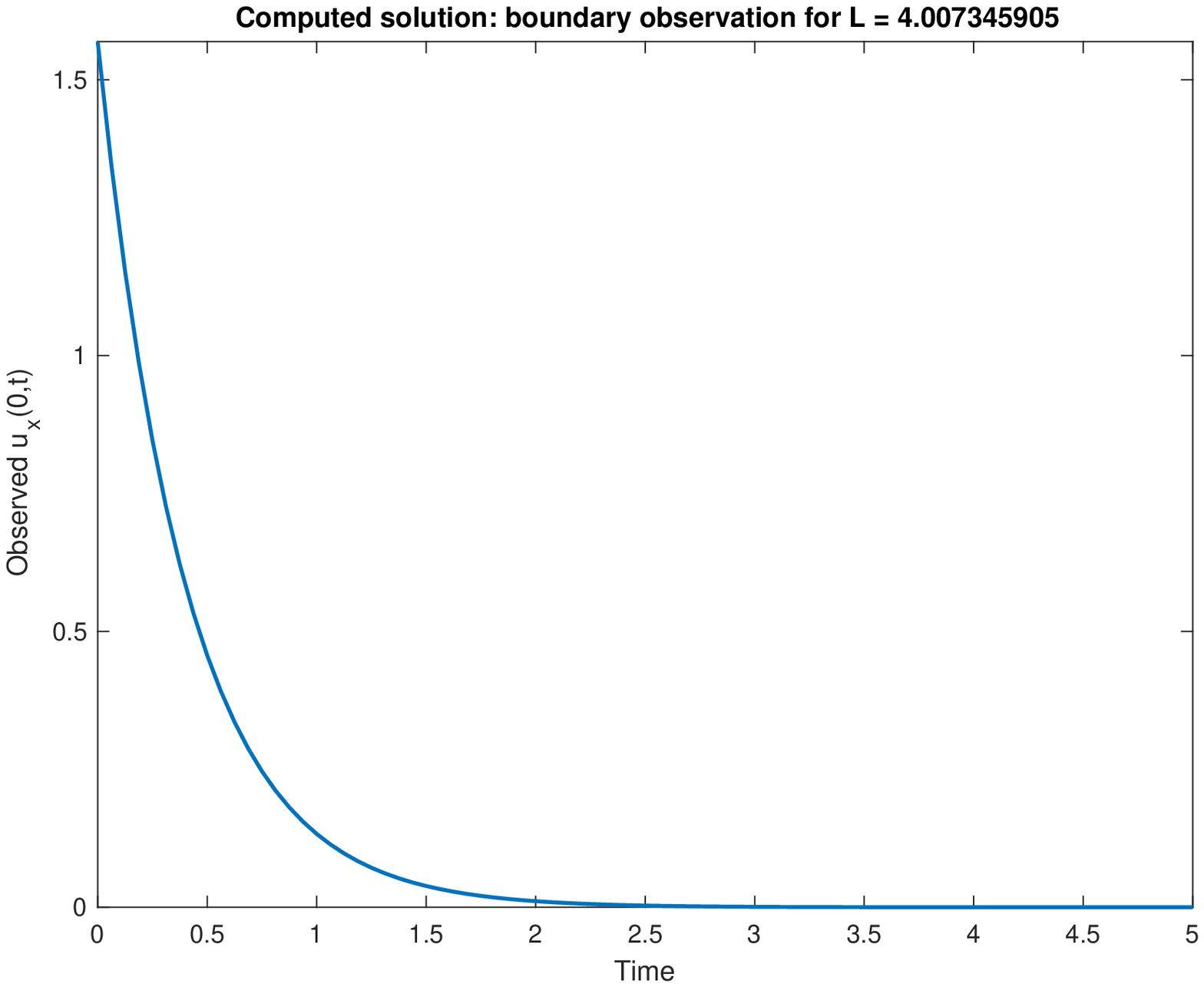}
\caption{Heat equation - fixed $u_0$ and $\eta=0$. The computed boundary observation $u_x(0,\cdot)$ for~$L_c^2=4.007345905$}
\label{Heat_Sol_case3b}
\end{minipage}
\end{figure}

   The corresponding computed boundary observations are displayed in~Figures~\ref{Heat_Sol_case3} and~\ref{Heat_Sol_case3b}, respectively.
   Thus, we confirm that these identical observations correspond, as we already knew, to different solutions. 
   
\subsection{Tests 2: The wave equation}\label{SSec-4.2}

   In this section, we consider the wave equation.
   The problem under study is the following:\\
   
   \textbf{Reformulation of IP-2:} \textit{Given $\eta = \eta(t)$, $u_0 = u_0(x)$, $u_1 = u_1(x)$, $T>0$ and $\beta = u_x(0,\cdot)$, find $\ell \in (\ell_0,\ell_1)$ such that 
   \begin{equation*}\label{pb.opt-w}
J(\ell) \le  J(\ell'),\quad \forall\, \ell'\in (\ell_0,\ell_1) ,
   \end{equation*}
where $J$ is given by~\eqref{eq.J}.
   Now, $u^\ell$ is the solution to~\eqref{problema-w} associated to the length~$\ell$.}
   
\

\noindent 
\textbf{Case 2.1: Wave equation with $(u_0,u_1) = (0,0)$ and  $\eta\neq 0$.} 

   We take $T=4$,  $\eta(t) =3\sin^3 t$ in~$(0,T)$, $u_0(x) \equiv 0$ and~$u_1(x) \equiv 0$.
   Starting from $L_{i} = 1.5$, we want to recover the desired length~$L_d= 2$.
   
\begin{figure}[h!]
\begin{minipage}[t]{0.49\linewidth}
\noindent
\includegraphics[width=\linewidth]{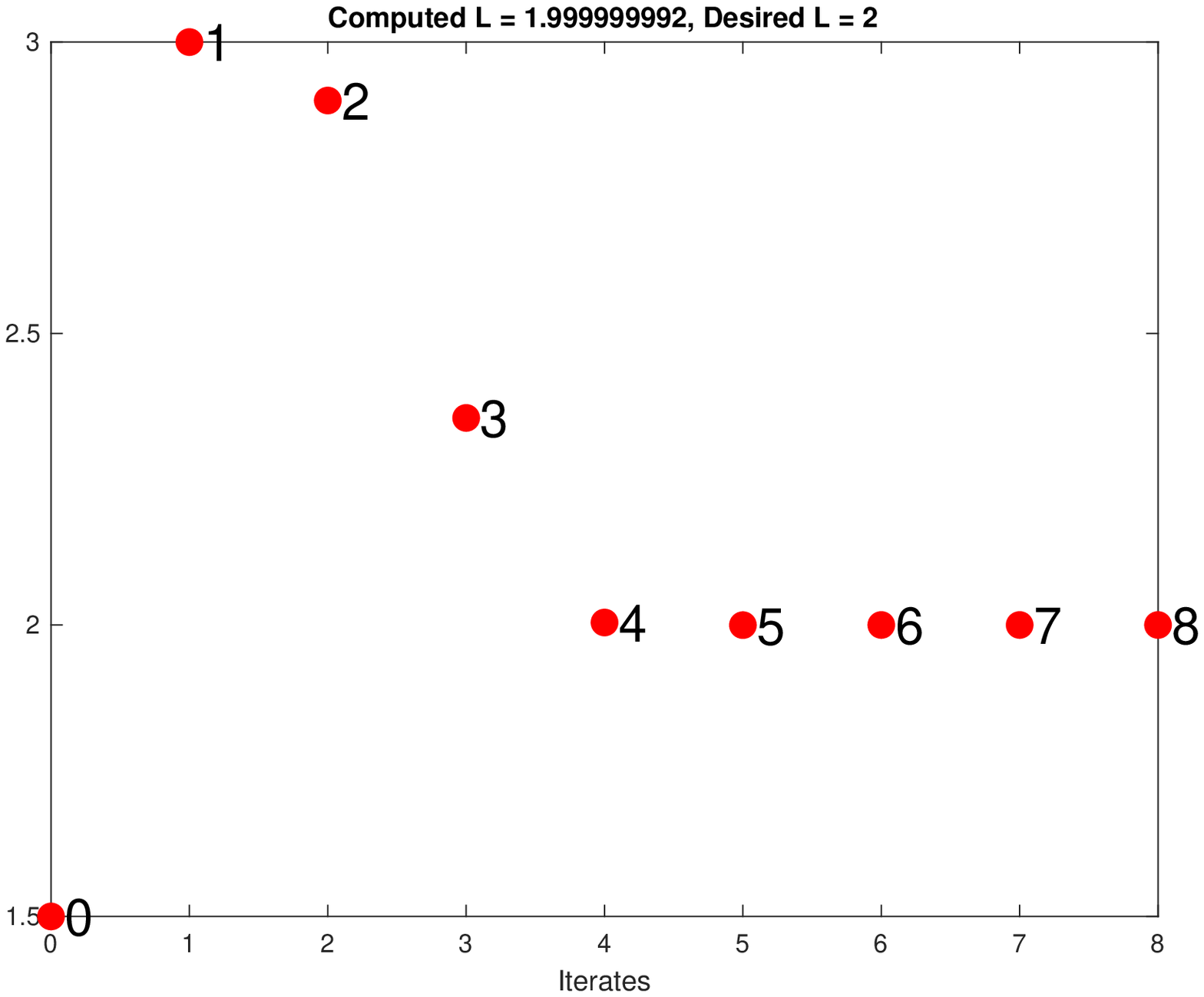}
\caption{Wave equation with~$(u_0,u_1) = (0,0)$ and~$\eta \not= 0$. The iterates in \texttt{active-set} algorithm.}
\label{Wave_iter_case1}
\end{minipage}
\hfill
\begin{minipage}[t]{0.49\linewidth}
\noindent
\includegraphics[width=\linewidth]{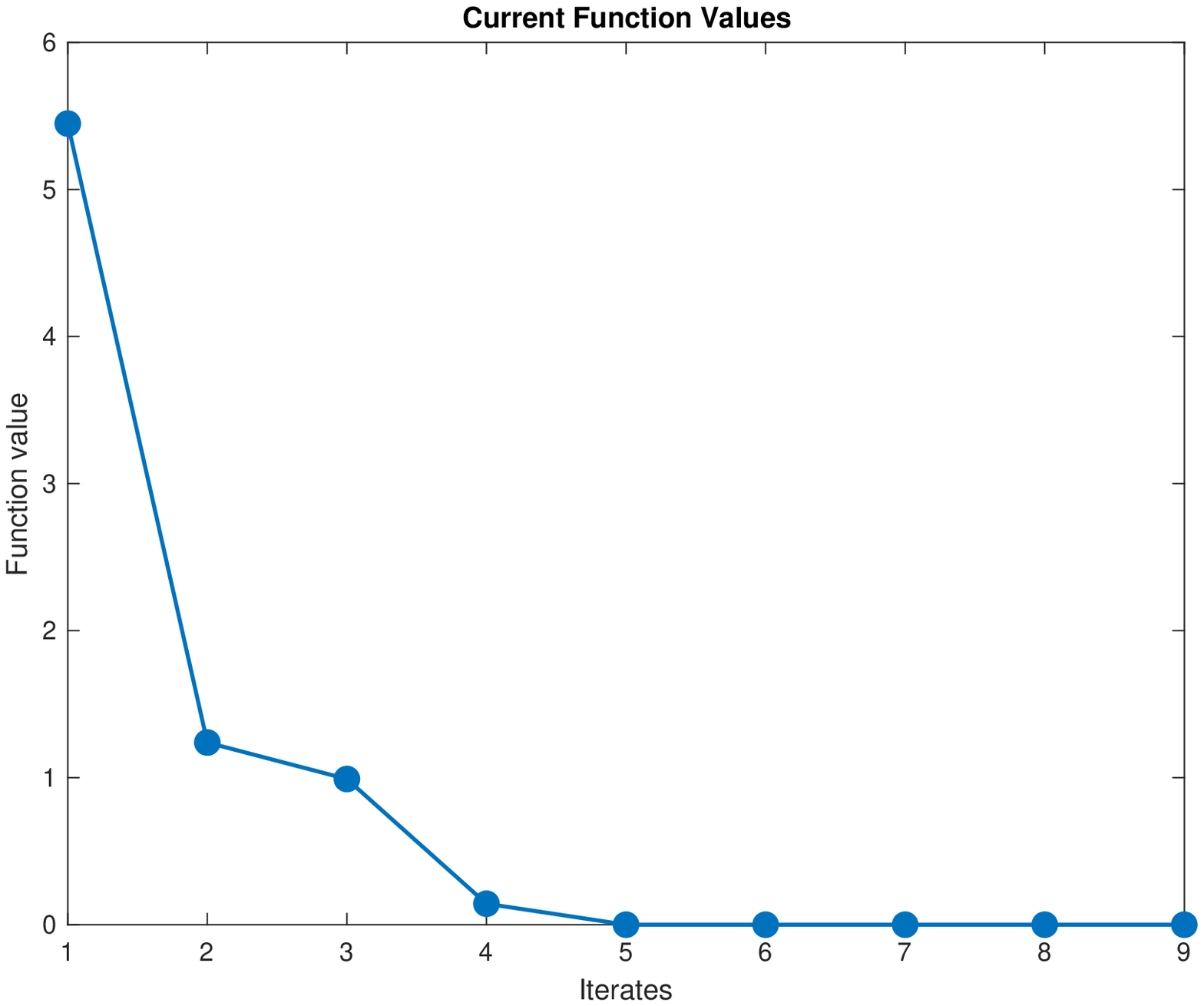}
\caption{Wave equation with~$(u_0,u_1) = (0,0)$ and~$\eta \not= 0$. Evolution of the cost.}
\label{Wave_J_case1}
\end{minipage}
\end{figure}

   Now, the numerical results are shown in Table~\ref{tab:wave_case1} and~Figures~\ref{Wave_Sol_case1}, \ref{Wave_iter_case1} and~\ref{Wave_J_case1}. 


%
\begin{table}[h!]
\begin{minipage}[t]{0.46\linewidth}
\centering
\renewcommand{\arraystretch}{1.2}
\medskip
\caption{Wave equation - $(u_0,u_1) = (0,0)$ and~$\eta \not= 0$. Results with random noise in the target
{(the desired length is $L_d = 2$). } \\
{ } \\
{ }}
\begin{tabular}{|cccc|} \hline
\%  noise  & Cost  & Iterates &  Computed $L_c$ \\
\hline
1\%           &   1.e-4        &       12       &   2.004287790  \\ \hline
0.1\%        &    1.e-6        &      8       &   2.000029532\\ \hline
0.01\%      &   1.e-8       &     8       &     2.000018464\\\hline
0.001\%    &   1.e-10      &     9         &    1.999995613\\ \hline
0\%           &   1.e-17       &     8         &    1.999999994  \\  \hline
\end{tabular}
\label{tab:wave_case1}
\end{minipage}
\hfill
\begin{minipage}[t]{0.54\linewidth}
\vspace{0cm}
\includegraphics[width=\linewidth]{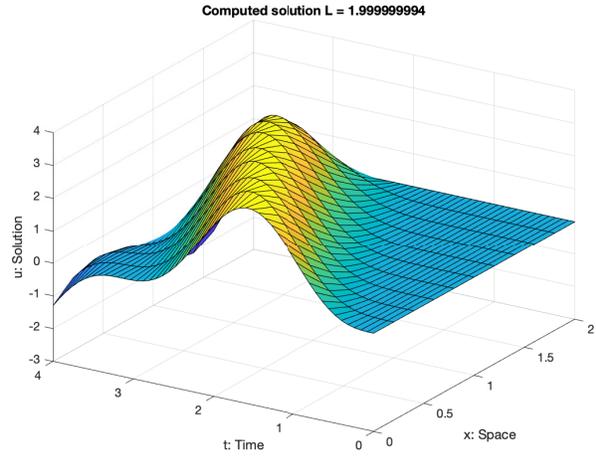}
\captionof{figure}{Wave equation with~$(u_0,u_1) = (0,0)$ and~$\eta \not= 0$. The computed solution.}
\label{Wave_Sol_case1}
\end{minipage}
\end{table}
%

\eject

\

\noindent 
\textbf{Case 2.2: Wave equation with $(u_0,u_1) \neq (0,0)$ and large $\eta$.}

   Here, $T=4$,  $\eta(t) = 3\sin^3 t$ in~$(0,T)$, $u_0(x) \equiv 0.4 \sin(\pi x)$ and~$u_1(x) \equiv 0$ and we start from~$L_{i} = 1.5$, trying to recover the desired length $L_d= 2$. 
   
\begin{figure}[h!]
\begin{minipage}[t]{0.49\linewidth}
\noindent
\includegraphics[width=\linewidth]{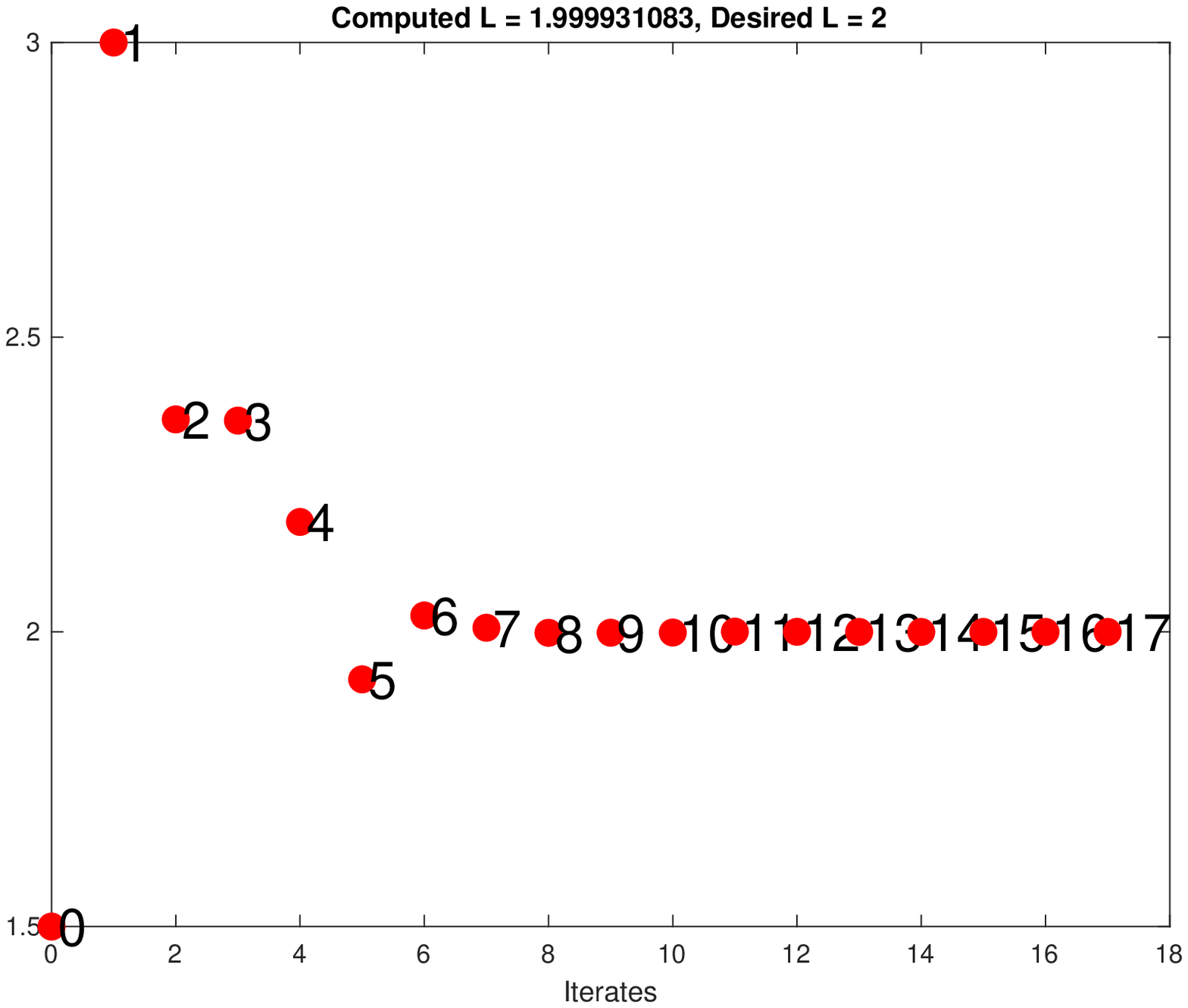}
\caption{Wave equation with~$(u_0,u_1) \not= (0,0)$ and~large~$\eta = 0$. The iterates in \texttt{active-set} algorithm.}
\label{Wave_iter_case2}
\end{minipage}
\hfill
\begin{minipage}[t]{0.49\linewidth}
\noindent
\includegraphics[width=\linewidth]{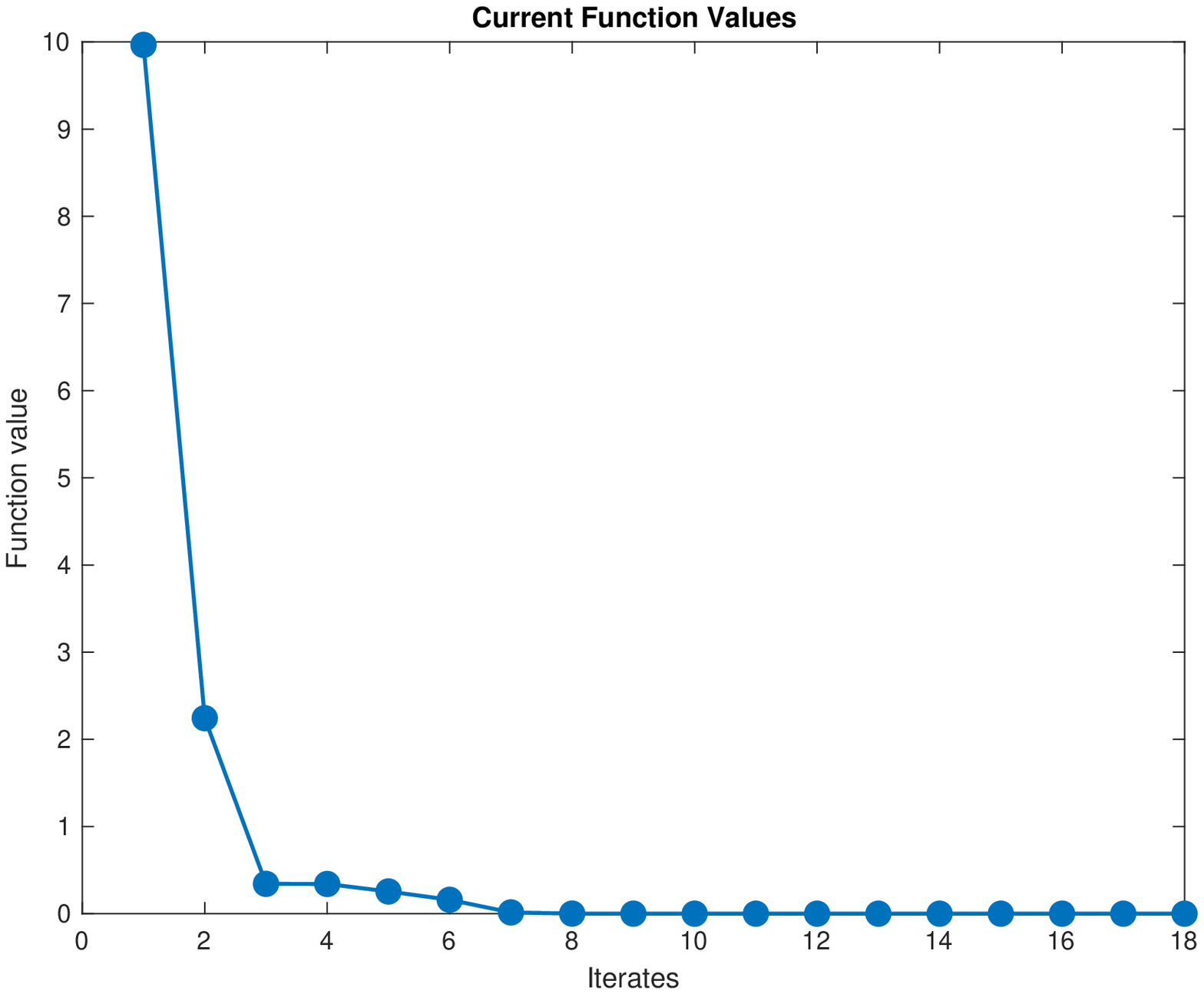}
\caption{Wave equation with~$(u_0,u_1) \not= (0,0)$ and~large~$\eta = 0$. The evolution of the cost.}
\label{Wave_J_case2}
\end{minipage}
\end{figure}

   See Table~\ref{tab:wave_case2} and~Figures~\ref{Wave_Sol_case2}, \ref{Wave_iter_case2} and~\ref{Wave_J_case2} 
for the numerical results.


%
\begin{table}[h!]
\begin{minipage}[t]{0.46\linewidth}
\centering
\renewcommand{\arraystretch}{1.2}
\medskip
\caption{Wave equation - $(u_0,u_1) \not= (0,0)$ and~large~$\eta = 0$. Cost and $L_c$ with random noise in the target
{(the desired length is $L_d = 2$). } \\
{ } \\
{ }}
\begin{tabular}{|cccc|} \hline
\%  noise  & Cost  & Iterates &  Computed $L_c$ \\
\hline
1\%           &   1.e-1        &       12       &   2.324511735  \\ \hline
0.1\%        &    1.e-6        &      11       &   2.000008724\\ \hline
0.01\%      &   1.e-7       &     9        &     1.999805367\\\hline
0.001\%    &   1.e-12      &     10         &    2.000000789\\ \hline
0\%           &   1.e-17       &     14         &    1.999931083  \\  \hline
\end{tabular}
\label{tab:wave_case2}
\end{minipage}
\hfill
\begin{minipage}[t]{0.54\linewidth}
\vspace{0cm}
\includegraphics[width=\linewidth]{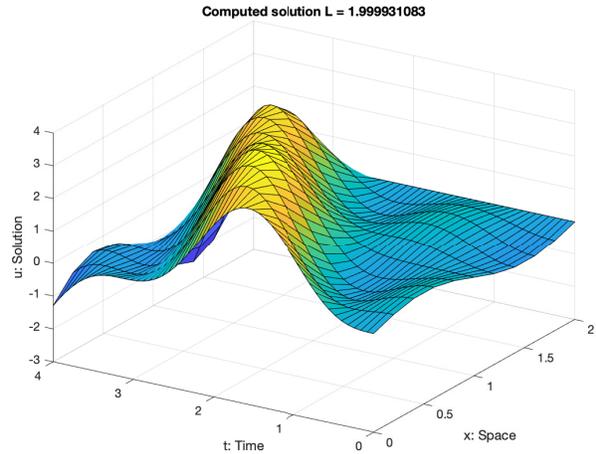}
\captionof{figure}{Wave equation with~$(u_0,u_1) \not= (0,0)$ and~large~$\eta = 0$. The computed solution.}
\label{Wave_Sol_case2}
\end{minipage}
\end{table}
%

\eject

\

\noindent 
\textbf{Case 2.3: Wave equation with $u_0\neq 0$ and ``small'' $\eta$.} 

   As in the previous section, we deal here with a case related to non-uniqueness. We take $T=4$ and the same data as in the Case 1.3.  
   
\begin{figure}[h!]
\begin{minipage}[t]{0.49\linewidth}
\noindent
\includegraphics[width=\linewidth]{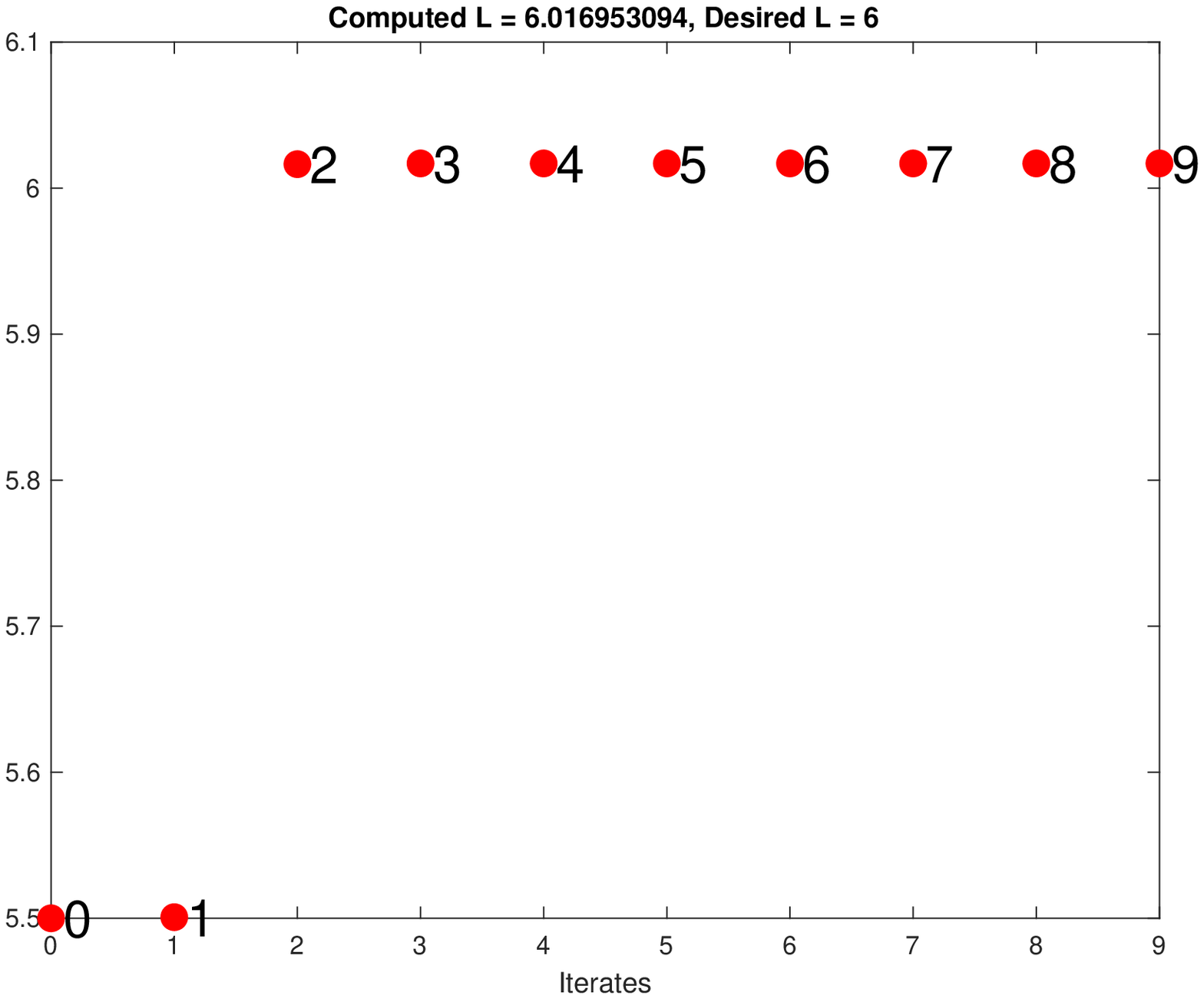}
\caption{Wave equation - $\eta=0$,  fixed $u_0(x)$. Iterates in \texttt{active-set} algorithm with $L^1_d=6$.}
\label{Wave_iter_case3}
\end{minipage}
\hfill
\begin{minipage}[t]{0.49\linewidth}
\noindent
\includegraphics[width=\linewidth]{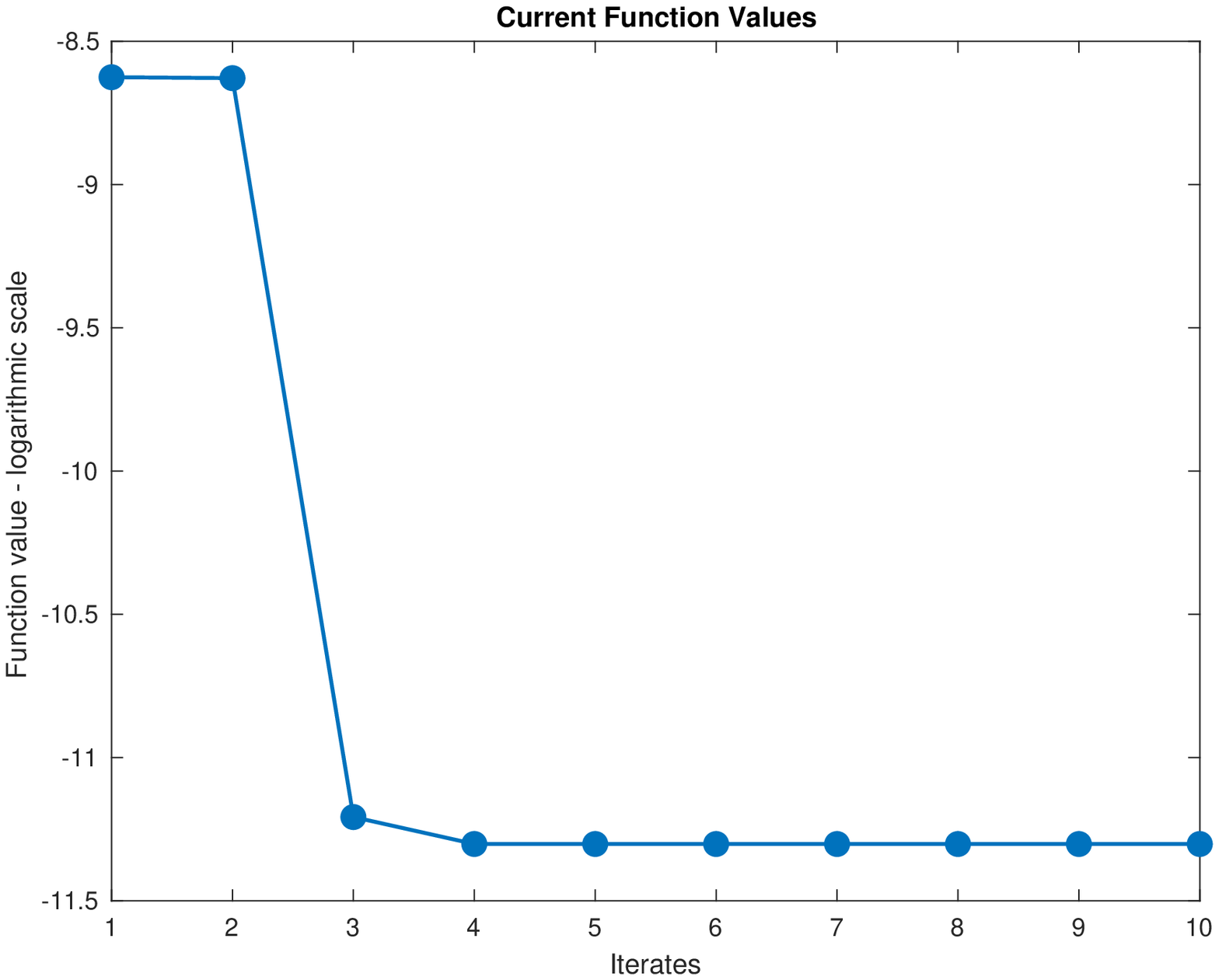}
\caption{Wave equation - $\eta=0$,  fixed $u_0(x)$. Evolution of the cost for $L^1_d=6$, $J(L^1_c)<10^{-6}$.}
\label{Wave_J_case3}
\end{minipage}
\end{figure}

\begin{figure}[h!]
\begin{minipage}[t]{0.49\linewidth}
\noindent
\includegraphics[width=\linewidth]{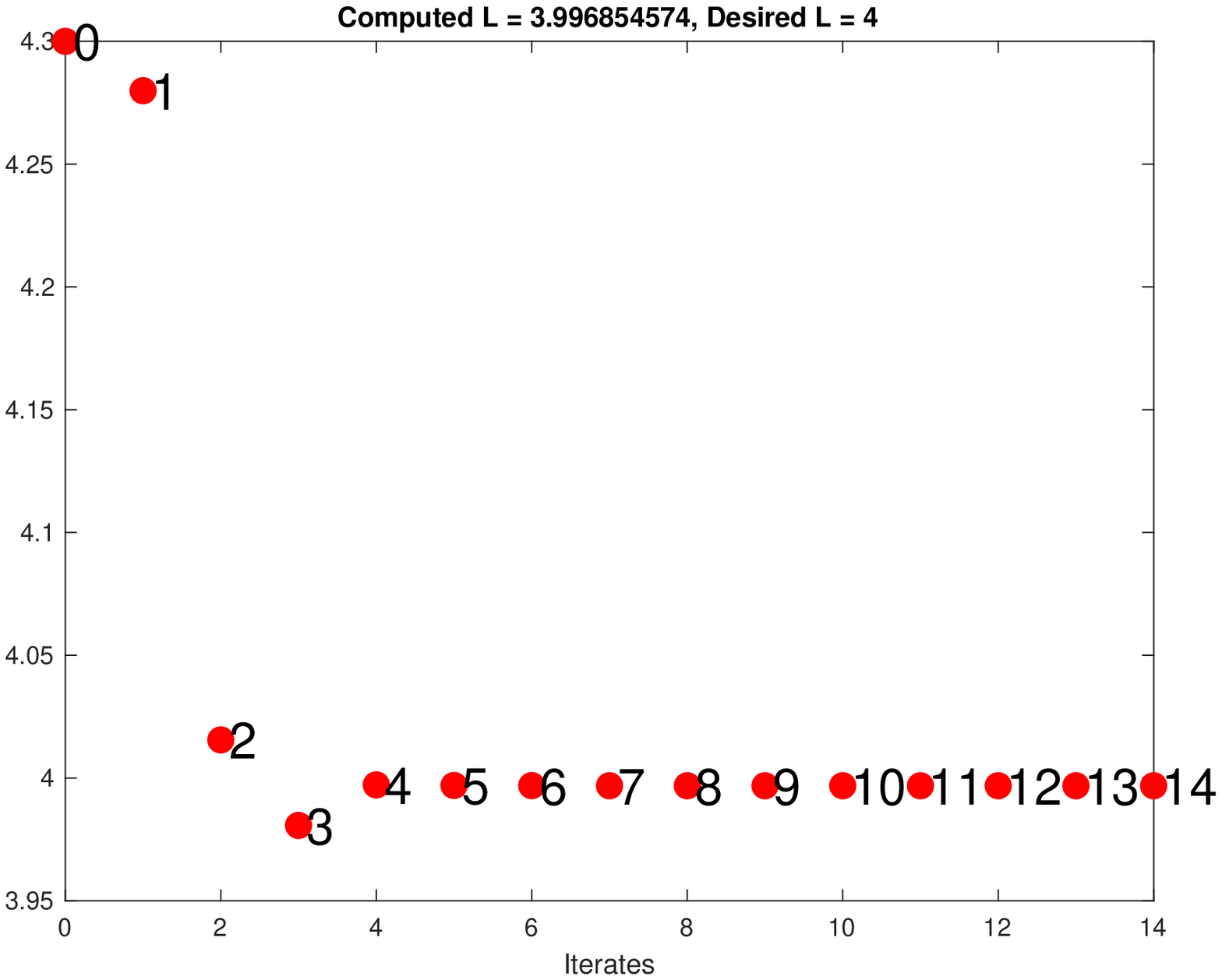}
\caption{Wave equation -  $\eta=0$,  fixed $u_0(x)$. Iterates in \texttt{active-set} algorithm with $L^2_d = 4$.}
\label{Wave_iter_case3b}
\end{minipage}
\hfill
\begin{minipage}[t]{0.49\linewidth}
\noindent
\includegraphics[width=\linewidth]{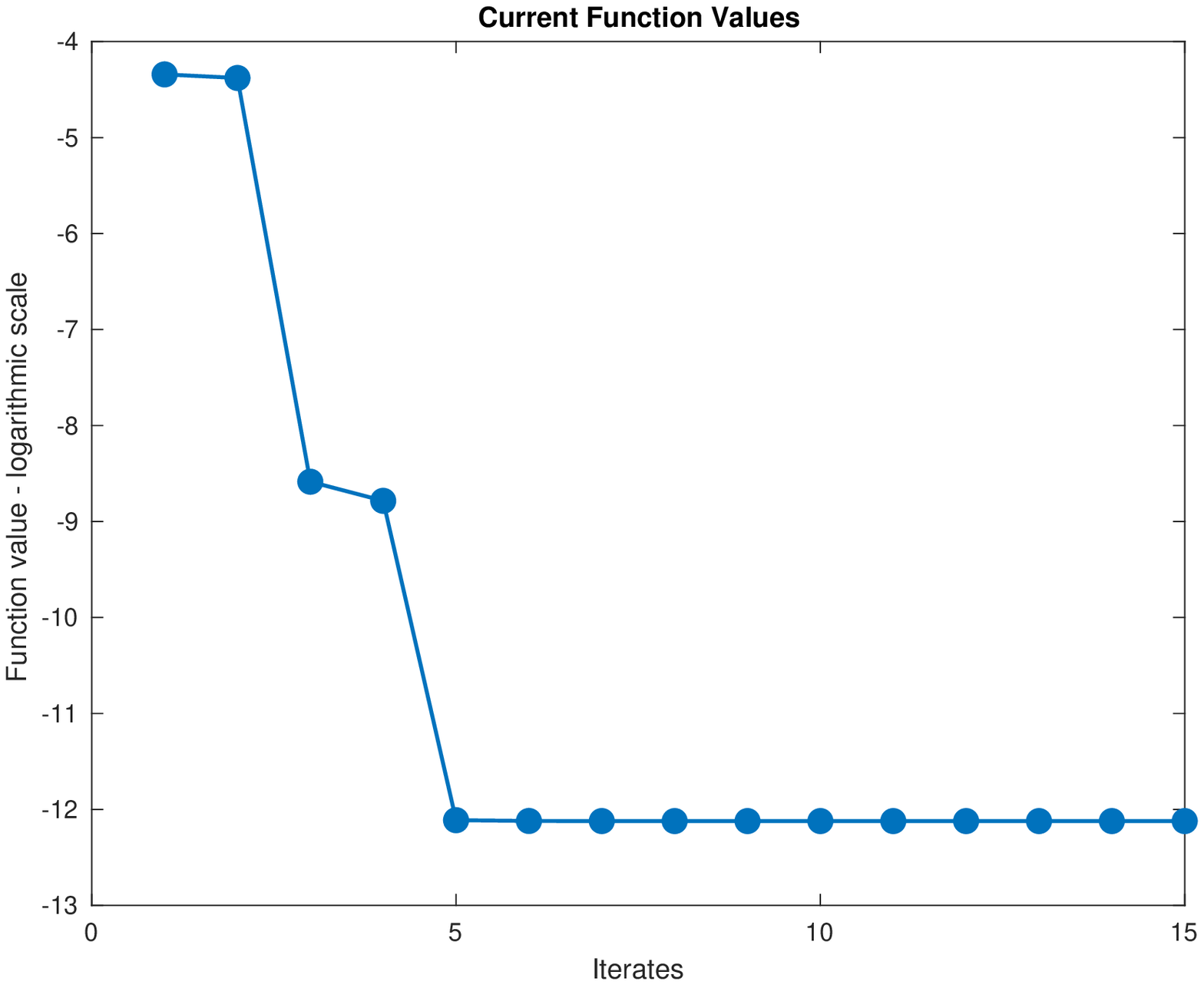}
\caption{Wave equation - $\eta=0$,  fixed $u_0(x)$. Evolution of the cost for $L^2_d=4$, $J(L_c^2)<10^{-6}$.}
\label{Wave_J_case3b}
\end{minipage}
\end{figure}

\begin{itemize}

\item Starting {\black from~$L_i = 5.5$,}  we get the results exhibited in~Figures~\ref{Wave_iter_case3} and~\ref{Wave_J_case3}, with~$L^1_c=6.016953094$ and associated final cost is $J(L_c^1)< 10^{-6}$.
   
\item Then, starting {\black from~$L_i = 4.5$,} we obtain the results exhibited in~Figures~\ref{Wave_iter_case3b} and~\ref{Wave_J_case3b}.
   The computed value is $L^1_c=3.996854574$ and the associated cost is again $J(L_c^1)< 10^{-6}$.
   
\end{itemize}

\begin{figure}[h!]
\begin{minipage}[t]{0.49\linewidth}
\noindent
\includegraphics[width=\linewidth]{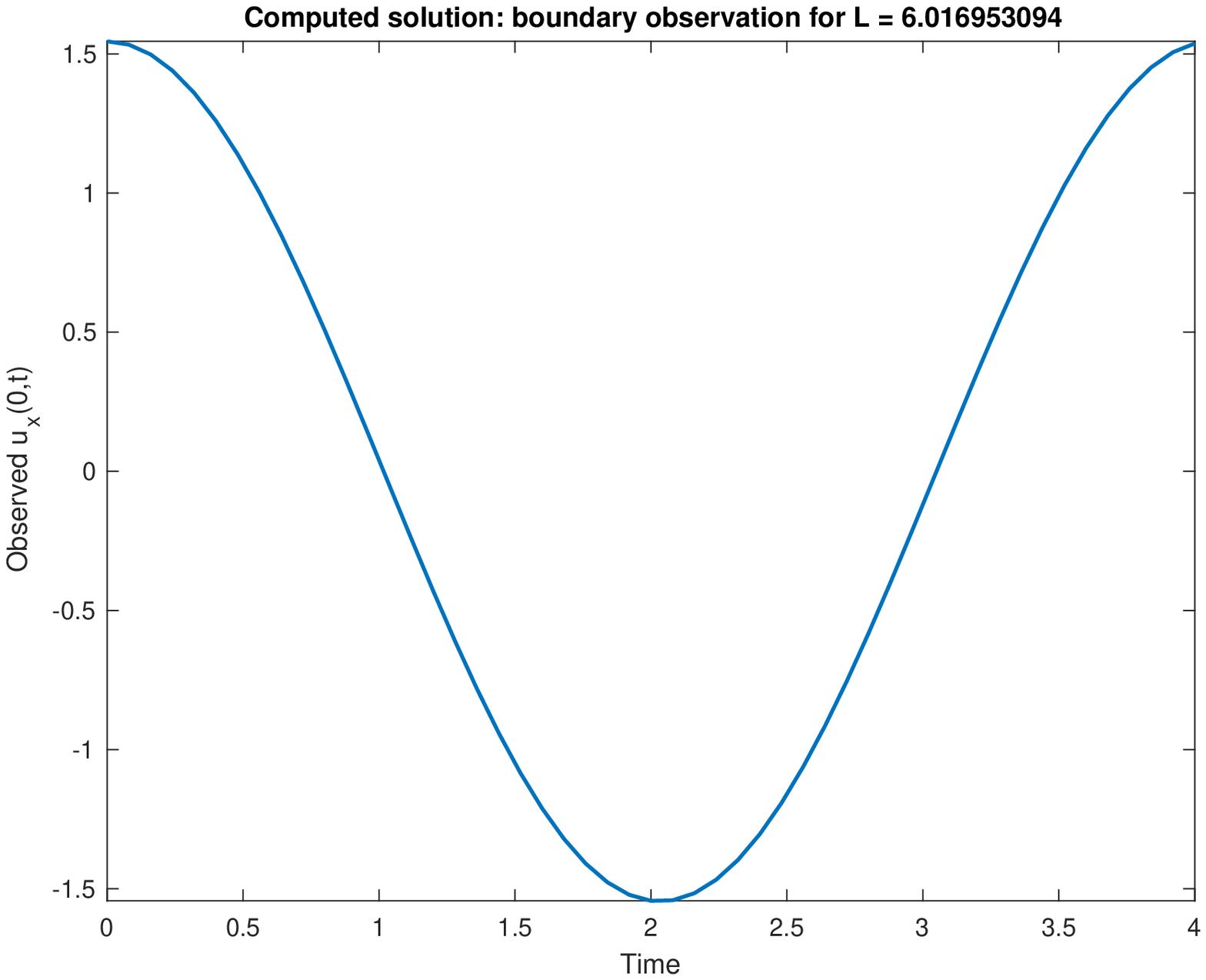}
\caption{Wave equation -  fixed $u_0(x)$, $\eta=0$. The computed boundary observation $u_x(0,t)$ corresponding to  $L^1_c=6.016953094$.}
\label{Wave_boObs_case3}
\end{minipage}
\hfill
\begin{minipage}[t]{0.49\linewidth}
\noindent
\includegraphics[width=\linewidth]{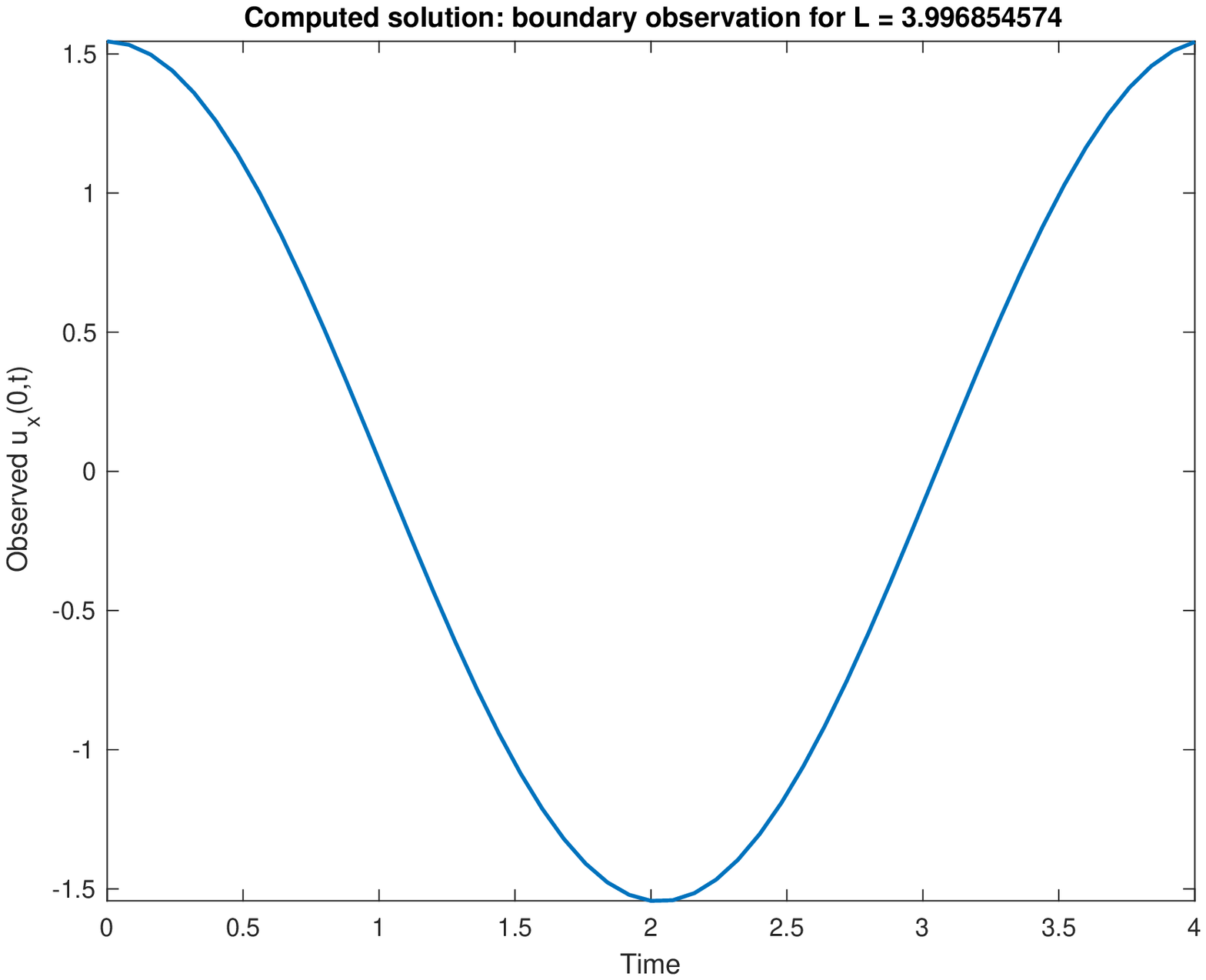}
\caption{Wave equation - fixed $u_0$, $\eta=0$. The computed boundary observation $u_x(0,t)$ corresponding to $L_c^2=3.996854574$}
\label{Wave_boObs_case3b}
\end{minipage}
\end{figure}

\begin{figure}[h!]
\begin{minipage}[t]{0.49\linewidth}
\noindent
\includegraphics[width=\linewidth]{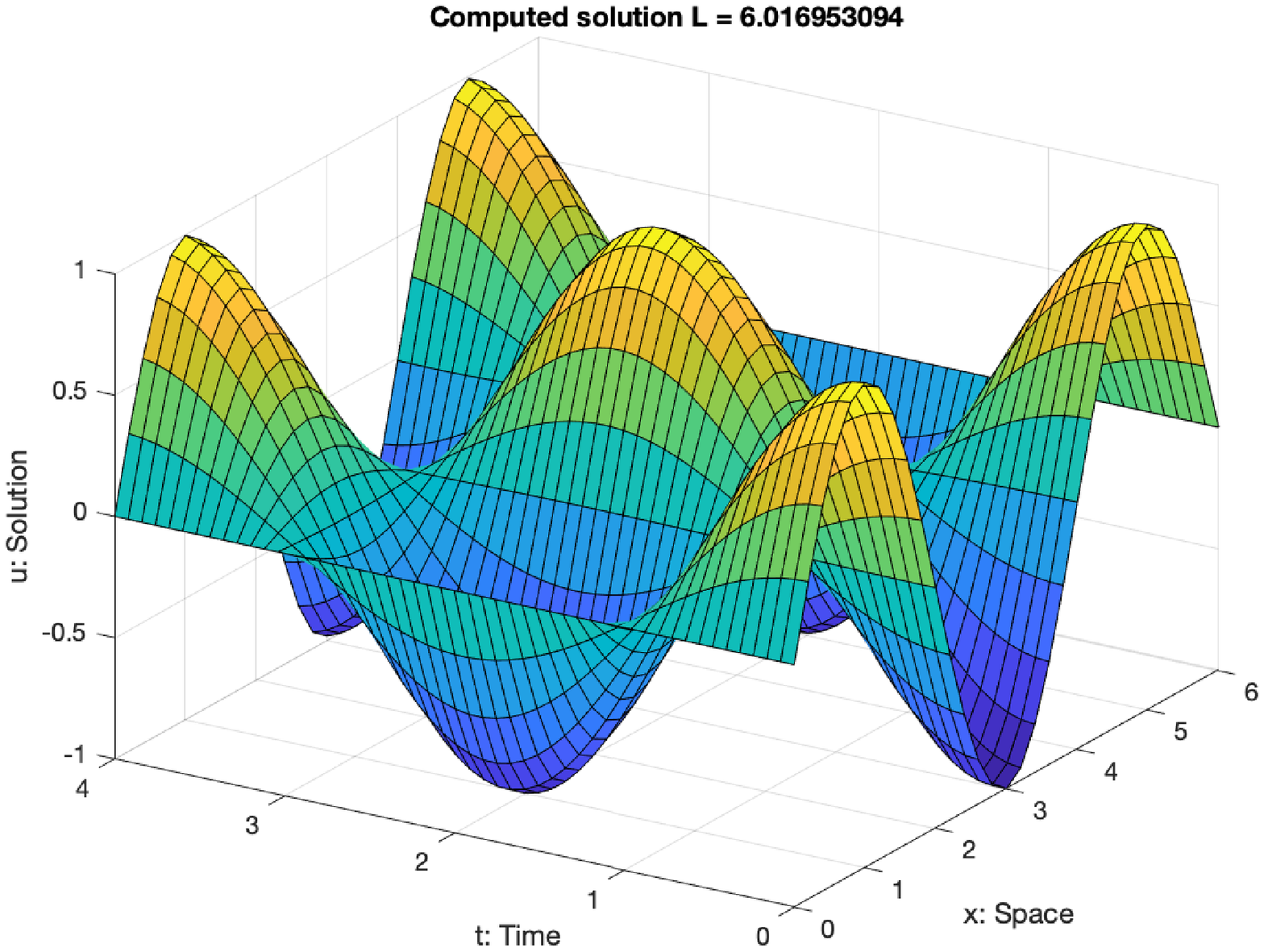}
\caption{Wave equation -  fixed $u_0(x)$, $\eta=0$. The computed solution corresponding to $L^1_c=6.016953094$.}
\label{Wave_Sol_case3}
\end{minipage}
\hfill
\begin{minipage}[t]{0.49\linewidth}
\noindent
\includegraphics[width=\linewidth]{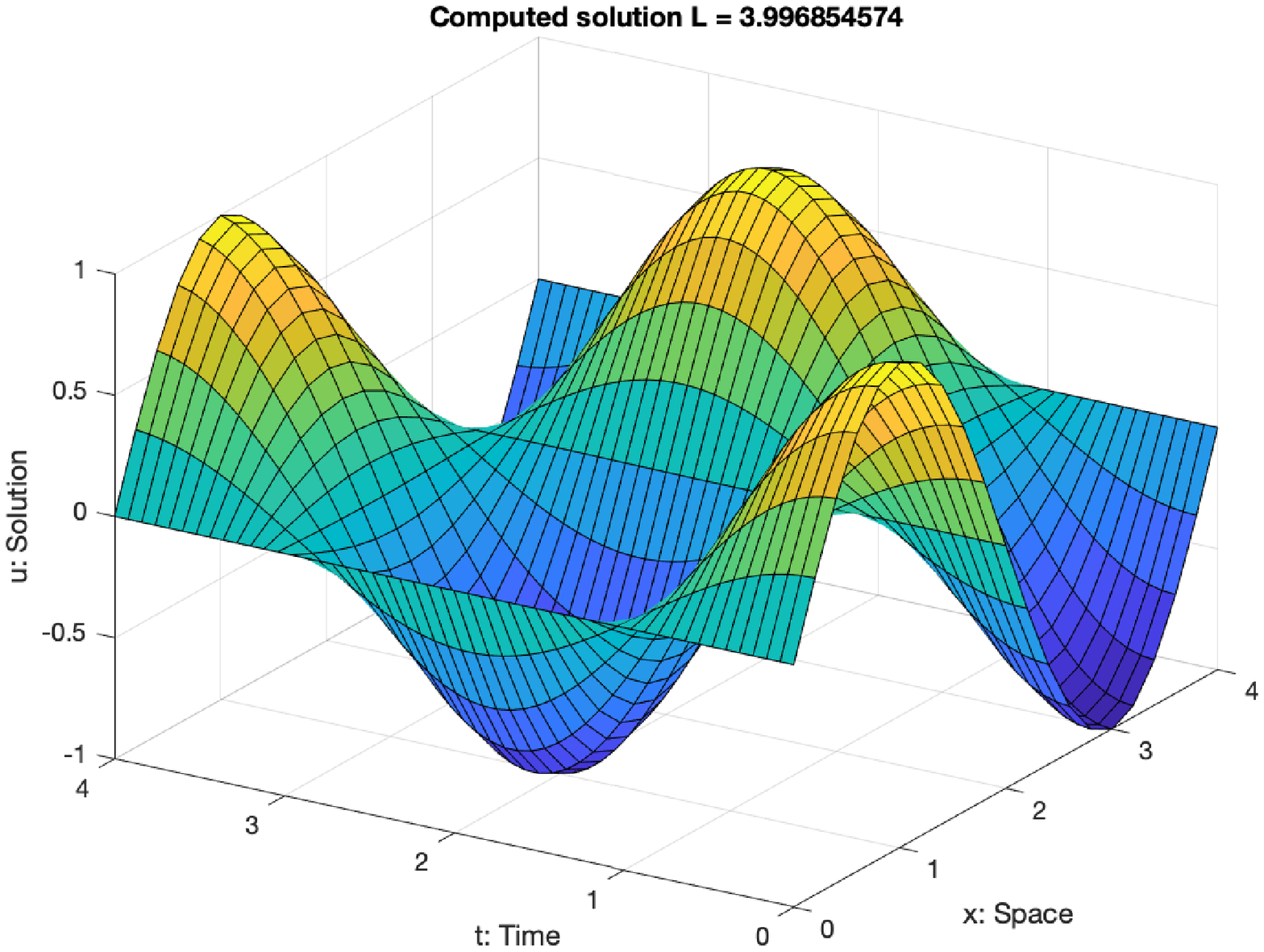}
\caption{Wave equation - fixed $u_0$, $\eta=0$. The computed solution corresponding to  $L_c^2=3.996854574$}
\label{Wave_Sol_case3b}
\end{minipage}
\end{figure}

   The computed solutions can be found in~Figures~\ref{Wave_Sol_case3} and~\ref{Wave_Sol_case3b}, respectively. In Figures~\ref{Wave_boObs_case3} and~\ref{Wave_boObs_case3b} we can appreciate boundary observations in the both cases, we see that the identical observations correspond to different solutions. 

\section{Some additional comments and open questions}\label{Sec-5}

   This section is devoted to complement the previous theoretical and numerical results.
   
   First, note that in the above inverse problems the observation has always been performed on the known boundary $x = 0$.
   It makes sense to work with information on the other lateral boundary of the domain.
   For instance, in the case of the heat equation the following problem is in order:
   
\

   \textbf{Problem~IP3 (heat equation, observation on the unknown boundary):} \textit{Fix $u_0 = u_0(x)$ and $\eta = \eta (t)$ in~\eqref{problema} in appropriate spaces and assume that $u_x|_{x=\ell}$ (the heat flow on the right) is known. Then, find $\ell$.}

\

   A similar formulation holds when the wave equation is considered.
   
   The results in Sections~\ref{Sec-2}--\ref{Sec-4} can be adapted to this new situation.
   For reasons of space, we leave the details to the reader.
   
   Other variants of the previous inverse problems can be obtained when different boundary conditions are imposed to the solution and/or different information is at hand.
   Thus, let us assume for example that the state satisfies
   \begin{equation*}\label{problema-N}
\begin{cases}
u_t - u_{xx} = 0,\, & 0 < x < \ell,\, 0 < t < T,\\
u_x(0,t)=\eta(t),\ \ u_x(\ell,t)=0,\, & 0 < t < T,\\
u(x,0)=u_0(x),\, & 0<x<\ell,
\end{cases}
   \end{equation*}
where $\eta$ and $u_0$ are given and let $u(0,\cdot)$ be given.
   Then, the related inverse problem is, once more, to recover~$\ell$.
   
   Again, arguing as above, theoretical uniqueness, non-uniqueness and stability results and numerical reconstruction results can be found.
   
   It is natural to try to extend the results in the previous sections to semilinear and nonlinear governing systems.
   For instance, it is significative to assume that $u$ solves the Burgers equation, together with the initial and boundary conditions in~\eqref{problema}.
   The analysis of the corresponding inverse problem will be carried out, together with other questions, in a forthcoming paper.
   
   Finally, let us mention that it makes sense to consider problems similar to those above in higher spatial dimensions.

\section*{Acknowledgements}

{\black
   The first author was supported by the Spanish Government's Ministry of Science, Innovation and Universities~(MICINN), under grant~PGC2018-094522-B-I00 and the Basque Government, under grant~IT12247-19.
   The third and fourth authors were partially supported by~MICINN, under grant~MTM2016-76690-P.

The fifth author was supported by Grant-in-Aid for Scientific Research (S)
15H05740 of Japan Society for the Promotion of Science and
by The National Natural Science Foundation of China
(no. 11771270, 91730303).
This work was prepared with the support of the "RUDN University Program 5-100".
}

\end{document}